\documentclass[12pt]{article}
 \usepackage{arydshln}
\renewcommand{\epsilon}{\varepsilon}

\usepackage{multicol}

\newcommand{\R}{\mathbb R}

\newcommand{\supp}{\mbox{supp}}
\newcommand{\diag}{\mbox{diag}}

\newcommand{\bea}{\begin{eqnarray}}
\newcommand{\eea}{\end{eqnarray}}
\newcommand{\be}{\begin{equation}}
\newcommand{\ee}{\end{equation}}
\newcommand{\beo}{\begin{equation*}}
\newcommand{\eeo}{\end{equation*}}
\newcommand{\beao}{\begin{eqnarray*}}
\newcommand{\eeao}{\end{eqnarray*}}
\newcommand{\ba}{\begin{array}}
\newcommand{\ea}{\end{array}}

\newtheorem{satz}{Theorem}[section]
\newtheorem{lem}[satz]{Lemma}

\newtheorem{kor}[satz]{Corollary}
\newtheorem{exam}[satz]{Example}

\usepackage[authoryear]{natbib}
\usepackage{amsmath}
\usepackage{amssymb}
\usepackage[ansinew]{inputenc}
\usepackage{graphicx}
\usepackage{epsfig}
\usepackage{xcolor}
\usepackage[draft]{changes}
\setremarkmarkup{(#2)}

\begin{document}

\parindent 0cm
\title{Optimal designs for  dose response
curves with common parameters}

\author{
{\small Chrystel Feller$^*$, Kirsten Schorning$^{\#}$, Holger Dette$^{\#}$}   \\
{\small $^*$Statistical Methodology } \\
{\small Novartis Pharma AG } \\
{\small 4002 Basel, Switzerland } \\
\and
{\small Georgina Bermann$^*$, Bj\"{o}rn Bornkamp$^*$ }   \\
{\small $^{\#}$Ruhr-Universit\"at Bochum } \\
{\small Fakult\"at f\"ur Mathematik } \\
{\small 44780 Bochum, Germany } \\
}

\maketitle
\begin{abstract}
A common problem in Phase II clinical trials is the comparison of dose response curves corresponding to different treatment groups. If the effect of the dose level is described by parametric regression models and the treatments differ in the administration frequency (but not in the sort of drug) a
reasonable  assumption is that the regression models for the different treatments share common parameters. \\
This paper develops optimal design theory for the comparison of different regression models with common parameters. We derive upper bounds on the number of support points of admissible designs, and explicit expressions for $D$-optimal designs are derived for 
frequently used dose response models with a common location parameter. If the location and scale parameter in the different models coincide, minimally supported designs are determined and sufficient conditions for their optimality in the class of all designs derived. The results are illustrated in a dose-finding study comparing monthly and weekly administration.
\end{abstract}


\section{Introduction} \label{sec1}
\def\theequation{1.\arabic{equation}}
\setcounter{equation}{0}

Adequately describing the dose-response relationship of a
pharmaceutical compound is of
paramount importance for achieving a successful clinical
development. \cite{sack:2014} recently conducted a review of the
reasons for delay or denial of approval of drugs by the Food and Drug
Administration (FDA). For those drug submissions that were not
approved in the first-time application, one of the most frequent
deficiencies was a statistical uncertainty related to the selected dose,
illustrating the importance of clearly determining an efficacious and
safe dose in Phase II dose-finding trials.

Efforts to improve this situation have  led to the
introduction of dose-response modeling approaches in a prospective
manner as the primary analysis method in dose finding studies, and
have become increasingly widespread in the past few years
[see among many others, \cite{grie:kram:2005}, \cite{bretz2005},
\cite{thom:2006}, \cite{drag:hsua:padm:2007}, \cite{phrma:2007},
\cite{thom:swee:soma:2014}]. These methods can more adequately address
the main questions of interest in Phase II dose-finding studies
(\textit{i.e.} determination of the dose-response curve and
estimation of target doses of interest) than AN(C)OVA based pairwise
comparisons. Moreover, it was pointed out by numerous authors that an
appropriate choice of the experimental conditions can improve
the statistical accuracy in dose-finding studies substantially. For this reason
there exists a large amount of literature discussing the problem of constructing optimal experimental designs
for regression models, which are commonly  used to describe the dose relationships
[see \cite{drag:hsua:padm:2007}, \cite{dett:bret:pepe:pinh:2008}, \cite{drag:fedo:wu:2008},
\cite{fangheda2008}, \cite{gilbert2010},
or \cite{bret:dett:pinh:2010}  among many others]. \\
 For many compounds a question closely related to ``dose'', the
amount of drug, is the administration frequency of the drug. In most
situations it is not adequate to assume that the same amount of drug
per time unit (\textit{e.g.} total daily dose) administered at different
dosing intervals
(\textit{e.g.} once a day or twice a day)
will lead to the same pharmacological effect. For example for once a
day administration the drug exposure inside the body will generally be
higher just after administration and lower just before the next
administration, compared to a twice a day administration, where the
same amount of drug is split into two doses in the morning and the
evening, leading to more uniform  drug exposure over the
day.

These considerations often lead to the need of
evaluating  the question of
finding the right dose as well as dosing frequency dose-finding
studies in Phase II.  One way of
modeling the dose-response curves in the different treatment groups
is to estimate the dose response curve corresponding to each of
them separately. This
can, however, be wasteful as certain aspects of the dose-response
curves for different group can be similar for both groups, suggesting a borrowing of
strength. When dose-response modeling is done in terms of parametric
dose-response models, one can often assume that certain
parameters of the dose-response curves for the two (or more) groups are shared, while other
parameters might be assumed to be  different between the curves. For example, if
the Emax function
\be \label{emax}
f (d, \theta_1, \theta^{(i)}_2)= \theta^{(i)}_0+ 
\frac {\theta^{(i)}_1 d}{\theta^{(i)}_2+d}, ~i=1,2
\ee
is used to model the dose response relationship for both groups
[see \cite{gabr:wein:2007} or \cite{thom:swee:soma:2014}], it is often
reasonable to assume that the placebo effect is the same between
groups, that is  $ \theta^{(1)}_0= \theta^{(2)}_0 =\vartheta_{11}$.
In some situations it might also make sense to assume that
the maximum efficacy for high doses is similar, i.e.
 $ \theta^{(1)}_1= \theta^{(2)}_1 =\vartheta_{12}$,
as a biological
maximum attainable effect might exist. However it might not be
adequate to assume that the dose providing half of the maximum
efficacy is the same for different treatment frequencies, which means $\theta^{(1)}_2  \neq \theta^{(2)}_2$.  The common parameters can then
be estimated more precisely allowing for a more accurate statistical analysis. An
example motivating the research of this paper can be found in Section  \ref{secEX2}.

The major question when planning such a dose-finding study then is
which doses to utilize in the different treatment groups and how to
split the total sample size between the groups. Statistically, this corresponds to the construction of optimal designs for different regression models (modeling the effect of the drug in the different groups) which share some common parameters.
To our best knowledge,  design problems in this case have not been considered
in the literature, and the goal of the  present paper is to derive optimal designs for such situations.
In Section \ref{secNOT} the model (in the context of $M$ treatment groups) is introduced  and the main
differences between the situation considered in the paper and the common optimal
design problems are explained. In Section \ref{secoptdes}  we derive some results on the comparison
of different designs for regression models with common parameters with respect to the Loewner ordering.
In particular we generalize recent results of admissible designs as presented in
\cite{yang2010}, \cite{detmel2011} and \cite{yangstuf2012b} and derive upper
bounds on the number of support points which cannot be improved upon in the Loewner ordering. Section \ref{secRES} is devoted to the construction of $D$-optimal designs
which are  well  suited for a ``global''  inference  as they  minimize the
maximum confidence interval length around the predicted dose-response
curve. Explicit expressions
for locally  $D$-optimal designs  for the commonly used dose response models are derived, if
some parameters of the models for the different groups coincide. We also discuss
minimally supported optimal designs and investigate if these designs are optimal within the
class of all designs. In Section \ref{secEX2} we illustrate the developed methods in a particular clinical
dose-finding study investigating two different treatment
groups. Finally, all technical details  and proofs are given
in Section \ref{sec:appendix} while Section \ref{sec:appendix2}  provides some more background
on the modeling problem discussed in Section \ref{secEX2}.
 \\
  For the sake of brevity and transparency, most parts of this paper consider   locally optimal designs
which require a-priori information about the unknown model parameters if the models are nonlinear
[see \cite{chernoff1953}]. In several situations  preliminary knowledge regarding the unknown parameters of a nonlinear model is available but not in a form that is accurate enough to specify one parameter guess. As illustrated in Section  \ref{secEX2}, locally optimal designs can
be used  as benchmarks for commonly used designs and also serve
as basis for constructing optimal
designs with respect to more sophisticated optimality
criteria,  which
are robust against a misspecification of the unknown parameters (and model) [see \cite{pronwalt1985}  or \cite{chaver1995}, \cite{dette1997} among others].
Following this line of research the methodology introduced in the present paper can be further developed to address uncertainty
in the preliminary information on the unknown parameters, and   we will illustrate this approach in Section \ref{secEX2},
where we also  discuss robust designs for the data example under consideration.

\section{Models with common parameters}\label{secNOT}
\def\theequation{2.\arabic{equation}}

Consider the regression models
 \begin{equation} \label{modgen}
 Y_{ij \ell} = f(d_j^{(i)}, \theta_1,  \theta^{(i)}_2) + \varepsilon_{ij \ell}  \quad i = 1, \dots, M; \, j= 1,\dots , k_i; \,  \ell = 1, \dots, n_{ij},
 \end{equation}
   where $\varepsilon_{ij \ell}$ are independent centered normally distributed random variables, i.e.\ $\varepsilon_{ij\ell} \sim \mathcal{N}(0,\sigma^2_i)$. 
 The assumption of a normal distribution in \eqref{modgen} is made for the sake of transparency. Other distributional
  assumptions can be treated exactly in the same way. 
   This means that $M$ different groups are investigated and in each group observations are taken at different experimental conditions $d^{(i)}_1, \ldots,   d^{(i)}_{k_i}$, which vary in possibly different design spaces, say  $\mathcal{X}_i=[0, d^{(i)}_{\max}] $ ($i = 1, \ldots, M$). At each dose level $d^{(i)}_j$ the experimenter can take $n_{ij}$ observations    and   $n_i= \sum_{j=1}^{k_i} n_{ij}$ denotes the number of observations in the $i$-th group ($i = 1, \ldots, M$).  Moreover, the total sample size is given by $n= \sum_{i=1}^M n_i$.
In general, the regression model $f(\cdot, \theta_1,  \theta^{(i)}_2)$ with a $(p+q)$-dimensional parameter vector $\theta^{(i)}=( \theta_1, \theta^{(i)}_2)$ is used to describe the dependency between the response and the effect in every group. We consider the same parametric form for all groups. Moreover,     the parameter vector $\theta_1 \in \R^p$ is assumed to be the same in all groups ($i=1, \ldots, M$), while $\theta^{(i)}_2 \in \mathbb{R}^q$ is different for different groups.
 Consequently, the  vector of unknown parameters is given by $\theta = (\theta_1, \theta^{(1)}_2, \ldots, \theta^{(M)}_2) \in \mathbb{R}^{m}$, where $m= p + qM$.
 The components of the vector are denoted by $\theta_1 =(\vartheta_{1}, \ldots , \vartheta_{p})$ and $ \theta^{(i)}_2 = ( \vartheta_1^{(i)},
 \ldots ,  \vartheta_q^{(i)})$ ($i=1,\ldots , M$). \\
Following  \cite{kiefer1974} we define for  $i=1, \ldots, M$  approximate designs $\xi_i$
 (on the design space $\mathcal{X}_i$) as probability measures with masses $\xi_{ij}$ at the experimental conditions $d^{(i)}_j \in \mathcal{X}_i  \ (j=1,\dots,k_i)$
  and a design $\mu$ as a probability measure on the set $\{ 1,\ldots, M \}$ assigning mass  $\lambda_i$ to the
 $i$th group. We collect these designs in the vector  $\xi = (\xi_1, \ldots, \xi_M, \mu)$, which is also called design (on the design
 space $ \mathcal{X}_1 \times \ldots \times  \mathcal{X}_M \times \{1,\ldots , M \}$) throughout
 this paper. If an approximate design  $\xi= (\xi_1, \ldots, \xi_M, \mu)$ is given and $N$ observations can be taken, a rounding procedure is applied to obtain integers $n_i$ and $n_{ij}$ ($i=1, \ldots M, \, j=1, \ldots, k_i$) from the not necessarily integer valued quantities $\lambda_i n$ and $\xi_{ij}n_i$, respectively  
 [see \cite{pukrie1992}].
 Then, under common assumptions of regularity and the assumption
\be \label{limit}
\lim_{n_i \rightarrow \infty} \tfrac{n_{ij}}{n_i}= \xi_{ij} \in (0, 1) \mbox{ and } \lim_{n \rightarrow \infty} \tfrac{n_i}{n} = \lambda_i \in (0, 1)
\ee
  $(i=1, \ldots, M$, $j= 1, \ldots, k_i)$,
 the maximum likelihood estimate
 $\hat \theta = ( \hat \theta_1, \hat \theta^{(1)}_2, \ldots, \hat \theta^{(M)}_2)$ satisfies (as $n \to \infty$)
 \beo
 \sqrt{n} (\hat \theta - \theta ) \xrightarrow{\mathcal{D}} \mathcal{N}( 0, M^{-1}(\xi, \theta)) ~,
 \eeo
 where the symbol $\xrightarrow{\mathcal{D}}$ denotes weak convergence. Here the matrix
\be \label{eq:info}
 M(\xi , \theta)= \int \int_{\mathcal{X}_z} h_z (d) h^T_z (d) d\xi_z(d) d\mu (z) =
 \sum^M_{i=1} \lambda_i M^{(i)} (\xi_i , \theta)
 \ee
 is called the information matrix of the design  $\xi = (\xi_1, \ldots, \xi_M, \mu)$ and will be derived in Section \ref{A1}. In \eqref{eq:info} the matrices $M^{(i)}$ are defined by
\be \label{eq:infoi}
M^{(i)} (\xi_i , \theta)=
  \int_{\mathcal{X}_i} h_i(d) h^T_i (d) d\xi_i(d)
  \ee
($i=1, \ldots , M$)  and
 \begin{eqnarray} \label{gradient}
h^T_i (d) = \frac{1}{\sigma_i} \Bigl( \tfrac{\partial}{\partial \theta_1} f(d, \theta_1, \theta^{(i)}_2),   \underbrace{0_q^T, \dots, 0_q^T}_{i-1} ,  \tfrac{\partial}{\partial  \theta^{(i)}_2}f(d, \theta_1, \theta^{(i)}_2), \underbrace{0_q^T, \dots, 0_q^T}_{M-i} \Bigr) \in \R^m
 \end{eqnarray}
 is the gradient of the function $f(d,\theta_1, \theta^{(i)}_2)$ with respect to the parameter $\theta \in \mathbb{R}^{m}$, where,
 $m=p+qM$, 
   $0_q \in \R^q$   denotes a vector with all entries equal to $0$.

\begin{exam}
{\rm
We assume that  $M=2$ and that the  regression functions $f(\cdot,\theta_1,  \theta^{(i)}_2)$ can be written as
\be  \label{locscal}
f(\cdot,\theta_1,  \theta^{(i)}_2)=  \vartheta_{1} + \vartheta_{2} f_0(\cdot, \theta^{(i)}_2)
 \ee
 with a given
   function $f_0$   [see \cite{bretz2005}]. Here   the location and scale parameters   $\theta_1 = (\vartheta_{1}, \vartheta_{2})^T \in \R^2$ are the same for all groups, while the parameters  $\theta^{(i)}_2 \in \R^q$ are  different. In this case we have $p=2$ and
 the vectors $h_1(d)$ and $h_1(d)$  are given by
\beao
 h^T_1(d) &=& \frac{1}{\sigma_1}
 \big(1,  f_0(d, \theta^{(1)}_2), \tfrac{\partial}{\partial \theta^{(1)}_2} f_0(d, \theta^{(1)}_2), 0_q^T \big), \\
 h^T_2(d) &=& \frac{1}{\sigma_2}
 \big(1,  f_0(d, \theta^{(2)}_2), 0_q^T , \tfrac{\partial}{\partial \theta^{(2)}_2} f_0(d, \theta^{(2)}_2)\big).
\eeao
As a further example, consider a   regression function  $f(\cdot,\theta_1,  \theta^{(i)}_2)$ of the form
\be  \label{eq:loc}
f(\cdot,\theta_1,  \theta^{(i)}_2)=  \theta_{1} + f_0(\cdot, \theta^{(i)}_2); \quad i=1,2,
 \ee
 with a given function $f_0$. If the location  parameter   $\theta_1$ is the same for the two groups and the parameters  $\theta^{(i)}_2 \in \R^q$ are  different,   we have $p=1$ and
 the vectors $h_1(d)$ and $h_1(d)$  are given by $ h^T_1(d) = \frac{1}{\sigma_1}(1, \tfrac{\partial}{\partial \theta^{(1)}_2} f_0(d, \theta^{(1)}_2), 0_q^T)$ and  $h^T_2(d)=  \frac{1}{\sigma_2}(1, 0_q^T , \tfrac{\partial}{\partial \theta^{(2)}_2} f_0(d, \theta^{(2)}_2))$.
 }
 \end{exam}

\section{Comparing designs in the Loewner ordering}  \label{secoptdes}
\def\theequation{3.\arabic{equation}}
\setcounter{equation}{0}

  \noindent
An optimal design  $\xi=(\xi_1, \ldots, \xi_M, \mu)$ maximizes a concave real valued function, say $\Phi$,  of the information matrix. Numerous
criteria have been proposed in the literature (see \cite{pukelsheim2006} among others) which can be used to discriminate
between competing designs and the particular case of $D$-optimality will be discussed in the subsequent section. The commonly used
optimality criteria are monotone with respect to the Loewner ordering, that is the relation $ M(\xi_1 , \theta)  \le M(\xi_2 , \theta)$ implies
$\Phi (M(\xi_1 , \theta) ) \le \Phi (M(\xi_2 , \theta) ) $.
For this reason we discuss at first some results for this ordering, which will be very  helpful for the explicit determination of
optimal designs in the following sections.

Throughout this paper let  $| {\cal A}|$ denote the cardinality of a set $ {\cal A}$ and we denote by $\supp(\xi)$ the support of the design $\xi=(\xi_1, \ldots, \xi_M, \mu)$. Moreover, we define the index $I(\xi_i)$ of the design $\xi_i$ on the interval $[0, d^{(i)}_{\max}]$ as the number of support points, where the boundary points $0$ and $d^{(i)}_{\max}$ are only counted by $1/2$ if they are support points of the design
$\xi_i$ ($i=1, \ldots, M$).\\
\noindent
 Note that the gradient \eqref{gradient} can be rewritten in the form
 \begin{equation} \label{pi}
 h_i(d) = \tfrac{1}{\sigma_i} \begin{pmatrix} I_{p \times p} & 0_{p \times q} \\
 												0_{(i-1)q \times p} & 0_{(i-1)q \times q} \\
												0_{q \times p} & I_{q \times q} \\
												0_{(M-i)q \times p} & 0_{(M-i)q \times p} \end{pmatrix} \begin{pmatrix}
\tfrac{\partial}{\partial \theta_1} f(d, \theta_1, \theta^{(i)}_2)\\
\tfrac{\partial}{\partial  \theta^{(i)}_2}f(d, \theta_1, \theta^{(i)}_2)
\end{pmatrix}
:= P_i \ g(d, \theta_1, \theta^{(i)}_2)
 \end{equation}
 where $P_i$ is a $(p+Mq)\times (p+q)$ block matrix, $ {I}_{p \times p}$  is the $p$-dimensional identity matrix and $g(d, \theta_1, \theta^{(i)}_2)$ is the $p+q$-dimensional gradient of $f(d, \theta_1, \theta^{(i)}_2)$ with respect to $(\theta_1,\theta^{(i)}_2)$ ($i=1, \ldots, M$).
 Consequently, for the information matrix \eqref{eq:infoi} the representation
 \begin{equation*}\label{eq:dlgi}
 M^{(i)}(\xi_i, \theta) = P_i \ \int_{\mathcal{X}_i} g(d, \theta_1, \theta^{(i)}_2) g(d, \theta_1, \theta^{(i)}_2) d\xi_i(d) \ P^T_i := P_i \ C(\xi_i, \theta_1, \theta^{(i)}_2) \ P^T_i
 \end{equation*}
holds,  where the  $(p+q) \times (p+q)$  matrix  $C(\xi_i, \theta_1, \theta^{(i)}_2)$ is defined by
\begin{eqnarray*}
C(\xi_i, \theta_1, \theta^{(i)}_2)&=&\int_{\mathcal{X}_i}
\begin{pmatrix}
 \Psi_{1,1}(d , \theta_1, \theta^{(i)}_2)& \ldots& \Psi_{1,p+q} (d, \theta_1, \theta^{(i)}_2)\\
 \vdots& \ddots & \vdots \\
\Psi_{p+q,1}(d, \theta_1, \theta^{(i)}_2) & \ldots & \Psi_{p+q,p+q}(d, \theta_1, \theta^{(i)}_2)
\end{pmatrix} d\xi_i(d)
\end{eqnarray*}
for $i=1, \ldots, M$.\\
In the following we will present a  generalization of results in
\cite{yang2010}, \cite{detmel2011} and \cite{yangstuf2012b}.
To be precise
 for $i=1, \ldots, M$ we define $\Psi_0(d) \equiv 1$ and choose a basis, say  $\{ \Psi_0(\cdot), \Psi^i_1(\cdot),  \ldots, \Psi^i_{k-1}(\cdot), \Psi^i_{k}(\cdot)\}$ for the space
 ${\rm span} (\Psi_{s,t}(\cdot, \theta_1, \theta^{(i)}_2) | 1 \leq s, t \leq p+q  \} \cup \{ 1\}  )$, where the dependence on the parameters is reflected by the upper index $i$ for the sake of a transparent notation. We also assume that  the function $\Psi^i_k(\cdot)$ is  a diagonal element of the matrix
 $C(\xi_i, \theta_1, \theta^{(i)}_2)$, does not coincide with any of the other elements $\Psi_{s,t}(\cdot, \theta_1, \theta^{(i)}_2) $
 and that $\{ \Psi_0(\cdot), \Psi^i_1(\cdot),  \ldots, \Psi^i_{k-1}(\cdot)\}$ is a basis of the space
$$
\mbox{span}\big(  \{ \Psi_{s,t} \mid s,t \in \{1,\ldots, p+q\}; \ \ \Psi_{s,t} \neq \Psi^{i}_k \} \cup \{1\}  \big) .
$$ 
For our first results we require the notation of Chebyshev system [see \cite{karstu1966}].  A set of $k$ real valued  functions $f_0,\ldots,f_{k-1}:[A,B] \rightarrow \mathbb{R}$ is called {\it Chebychev system} on the interval $[A,B]$ if and only if it fulfills the inequality
$$\det
\begin{pmatrix}
f_0(x_0)&\ldots &f _0(x_{k-1})\\
\vdots & \ddots & \vdots \\
f_{k-1}(x_0)&\ldots &f_{k-1}(x_{k-1})
\end{pmatrix}
>0
$$
for any points $x_0,\ldots,x_{k-1}$ with $A \leq x_0<x_1\ldots<x_{k-1} \leq B$.

\begin{lem}\label{theo:garza}   ~\\
(1) If  for all $i=1,\ldots,M$ the sets $\{\Psi_0(\cdot),  \Psi^i_1(\cdot), \ldots, \Psi^i_{k-1}(\cdot)\}$  and 
$\{ \Psi_0(\cdot), \Psi^i_1(\cdot), \ldots, \Psi^i_{k-1}(\cdot), \Psi^i_{k}(\cdot)\}$ are Chebychev systems  on the interval $\mathcal{X}_i =[0,d^{(i)}_{\max}]$,
then for any design $\xi$ there exists a design $\xi^+= (\xi^+_1, \ldots, \xi^+_M, \mu)$ with $|\supp(\xi_i^+)|\leq \tfrac{k+2}{2} \ (i=1,\ldots,M)$, such that $M(\xi^+, \theta) \geq M(\xi, \theta)$. If the index of the design $\xi_i$ satisfies $I(\xi_i) < \tfrac{k}{2}$ the design coincides with the design $\xi$. In the case $I(\xi_i) \geq \tfrac{k}{2}$, the following two assertions are valid. 
\begin{enumerate}
\item[(1a)] If $k$ is odd, then $\xi^+_i$ has at most $\tfrac{k+1}{2}$ support points and $\xi^+_i$ can be chosen such that its support contains $d^{(i)}_{\max}$ ($i=1, \ldots, M$).
\item[(1b)] If $k$ is even, then $\xi^+_i$ has at most $\tfrac{k+2}{2}$ support points and $\xi^+_i$ can be chosen such that its support contains the points $0$ and $d^{(i)}_{\max}$ ($i=1, \ldots, M$).
\end{enumerate}
(2)   If  for all $i=1,\ldots,M$ the sets $\{ \Psi_0(\cdot), \Psi^i_1(\cdot), \ldots, \Psi^i_{k-1}(\cdot)\}$  and
$\{  \Psi_0(\cdot), \Psi^i_1(\cdot), \ldots, \Psi^i_{k-1}(\cdot), -\Psi^i_{k}(\cdot)\}$ are Chebychev systems  on the interval $\mathcal{X}_i =[0,d^{(i)}_{\max}]$,
then for any design $\xi$ there exists a design $\xi^-= (\xi^-_1, \ldots, \xi^-_M, \mu)$ with
$|\supp(\xi_i^-)|\leq \tfrac{k+2}{2}$ $(i=1,\ldots,M)$, such that $M(\xi^-, \theta) \geq M(\xi, \theta)$.  If the index of the design $\xi_i$ satisfies $I(\xi_i) < \tfrac{k}{2}$ the design coincides with the design $\xi$. In the case $I(\xi_i) \geq \tfrac{k}{2}$, the following two assertions are valid. 
\begin{enumerate}
\item[(2a)] If $k$ is odd, then $\xi^-_i$ has at most $\tfrac{k+1}{2}$ support points and $\xi^-_i$ can be chosen such that its support contains $0$.
\item[(2b)] If $k$ is even, then $\xi^-_i$ has at most $\tfrac{k}{2}$ support points.
\end{enumerate}
\end{lem}

Lemma \ref{theo:garza} provides an upper bound for the maximal number of support points if the functions $\Psi_0(\cdot), \Psi^i_1(\cdot), \ldots, \Psi^i_{k-1}(\cdot)$  and
$\Psi_0(\cdot), \Psi^i_1(\cdot), \ldots, \Psi^i_{k-1}(\cdot), \Psi^i_{k}(\cdot)$ are Chebychev systems for the different   groups $i=1, \ldots, M$. Note that this bound is the same independently from the dimension of $\theta_1$, since the number of support points is bounded  in every group $1, \ldots, M$ separately.  
The next lemma  shows that (for the commonly used dose response models) it is sufficient to allocate
only patients from the group with the smallest population variance to placebo.

\begin{lem}\label{lem:placebo}
Assume that the design spaces are given by
$\mathcal{X}_i =[0,d^{(i)}_{\max}]$ $(i=1,\ldots , M)$ and that
 the regression models  are  given by \eqref{locscal} or by \eqref{eq:loc},
 where the function   $f_0$ is differentiable with respect to $\theta^{(i)}_2$ $(i=1,\ldots , M)$. Moreover, assume that $f_0(0, \theta_2) = 0$ and $\tfrac{\partial}{\partial  \theta_2} f_0 (0,  \theta_2) = 0$.
  If $\eta = (\eta_1,\ldots , \eta_M, \nu)$ denotes a design with $0 \in \mbox{supp} (\eta_j)$   for  (at least) one index $j$, then there
 exists a design $\xi = (\xi_1,\ldots , \xi_M, \mu)$ with the following properties
 \begin{eqnarray*}
M(\eta , \theta)  \le M(\xi , \theta), ~ 0 \in \mbox{supp} (\xi_{j^*}), ~ 0 \not \in \mbox{supp} (\xi_{j}) ~~\mbox{for all } j  \not = j^*,
\end{eqnarray*}
where
 $j^* \in \mbox{argmin}_{i=1, \ldots, M} \sigma^2_i$.
\end{lem}
In the following discussion we will apply   the previous results to some of the commonly used  dose response models,
namely the Emax model, linear-in-log and exponential model
 [see \cite{gabr:wein:2007}], which are  listed in Table  \ref{tab:consmod}.
 In this table we also illustrate our notation again. The left part of the table  corresponds to a model with a common location parameter
 (namely $\theta_1$),
 while the right part of the table shows a model with a
 common location ($\vartheta_{1} $) and scale  parameter ($\vartheta_{2}$).
We note that all these models satisfy  the conditions of Lemma \ref{theo:garza} and Lemma \ref{lem:placebo}.


\begin{table}[t]
\centering
\fbox{
\begin{minipage}{0.2\linewidth}
\begin{equation*}
\mbox{model}
\end{equation*}
\begin{equation*}
\mbox{Emax}
\end{equation*}
\begin{equation*}
\mbox{Linear-in-log}
\end{equation*}
\begin{equation*}
\mbox{Exponential}
\end{equation*}
\end{minipage}
\begin{minipage}{0.32\linewidth}
\begin{equation*}
\mbox{location}
\end{equation*}
\begin{equation}\label{eq:Emax2}
\theta_{1} + \vartheta^{(i)}_{1} \tfrac {  d}{\vartheta^{(i)}_2+d}
\end{equation}
\begin{equation} \label{eq:log2} 
\theta_{1} +\vartheta^{(i)}_{1} \log ( {d/  \vartheta^{(i)}_2} +1 )
\end{equation}
\begin{equation}\label{eq:exp2} 
 \theta_{1} +\vartheta^{(i)}_{1} (\exp (  {d} / {\vartheta^{(i)}_{2}} ) -1)
\end{equation}
\end{minipage}
\begin{minipage}{0.35\linewidth}
\begin{equation*}
\mbox{location and scale}
\end{equation*}
\begin{equation}\label{eq:Emax2a} 
\vartheta_{1} + \vartheta_{2} \tfrac {d}{\theta^{(i)}_2+d}
\end{equation}
\begin{equation} \label{eq:log2a} 
 \vartheta_{1} +\vartheta_{2} \log ( {d/  \theta^{(i)}_2} +1 )
\end{equation}
\begin{equation}\label{eq:exp2a} 
 \vartheta_{1} + \vartheta_{2} \Bigl(\exp({d}/{\theta^{(i)}_2}) -1 \Bigr) 
\end{equation}
\end{minipage}}
\caption{\it Commonly used dose response models for $i=1, \ldots, M$. Left column: The placebo effect is the same in every group
(common location). Right column: Both the placebo effect and the scale parameter coincide in every group (common location and scale).}
\label{tab:consmod}
\end{table}

\begin{kor}\label{kor:garza}
Let  $\xi=(\xi_1, \ldots, \xi_M, \mu)$ denote an arbitrary design with  $|\supp(\xi_i)|\geq 3$ $(i=1,\ldots,M)$ and
assume (w.l.o.g) that $\sigma^2_1= \min_{i=1, \ldots, M} \sigma^2_i$.
\begin{itemize}
\item[(1)] If the regression model is given by the Emax model \eqref{eq:Emax2}  or  \eqref{eq:Emax2a}, then  there exists a design $\xi^+=(\xi^+_1, \ldots, \xi^+_M, \mu)$
with at most $2M +1 $ support points such that $M(\xi^+, \theta) \geq~M(\xi, \theta)$.  Moreover,  $\xi^+$ can be chosen such that  $|\supp(\xi_1^+)|= 3$ with  $0,~d^{(1)}_{\max} \in \supp(\xi_1^+)$ and $|\supp(\xi_i^+)|=2$ with
 $d^{(i)}_{\max} \in \supp(\xi^+_i)$ ($i=2, \ldots, M$).
\item[(2)] If the regression model is given by the linear-in-log model  \eqref{eq:log2}  or  \eqref{eq:log2a}, then there exists a design $\xi^+=(\xi^+_1, \ldots, \xi^+_M, \mu)$ with at most $2M +1 $ support points such that $M(\xi^+, \theta) \geq M(\xi, \theta)$. Moreover,  $\xi^+$ can be chosen such that  $|\supp(\xi_1^+)|= 3$ with  $0,~d^{(1)}_{\max} \in \supp(\xi_1^+)$ and $|\supp(\xi_i^+)|=~2$ with
  $d^{(i)}_{\max} \in \supp(\xi^+_i)$ ($i=2, \ldots, M$).
\item[(3)] If the regression model is given by the exponential model \eqref{eq:exp2} or \eqref{eq:exp2a}, then there exists a design
$\xi^+=(\xi^+_1, \ldots, \xi^+_M, \mu)$ with at most $3M$ support points such that $M(\xi^+, \theta) \geq M(\xi, \theta)$. Moreover, $\xi^+_i$ can be chosen such that $|\supp(\xi_i^+)|= 3$  and
 $d^{(i)}_{\max} \in \supp(\xi^+_i)$ $(i=1,\ldots,M)$.
\end{itemize}
\end{kor}

\section{$D$-optimal designs}\label{secRES}
\def\theequation{4.\arabic{equation}}
\setcounter{equation}{0}

When one of the major
purposes of the study is to determine the dose-response curve,
$D$-optimal designs are well suited  as they  minimize the
maximum confidence interval length around the predicted dose-response
curve [see \cite{silvey1980}].
 Following \cite{chernoff1953},  a design
 $\xi = (\xi_1, \ldots, \xi_M, \mu)$ is called (locally) $D$-optimal for the information matrix given in \eqref{eq:info}
  if it maximizes the determinant of the information matrix 
   $\det (M(\xi, \theta))$ in the class of all designs  $\xi $ on $\mathcal{X}_1 \times \ldots \times \mathcal{X}_M \times \{1, \ldots, M\}$.
A main tool of optimal design theory are equivalence theorems which, on the one hand provide a
simple checking condition for the optimality of a given  design, and on the other hand, are the basis
of many procedures for their numerical construction.
 Moreover, these characterizations of optimality can also be used to derive structural properties of optimal designs.
 The following result provides the equivalence theorem for the
$D$-optimality criterion corresponding to the matrix given in \eqref{eq:info}.
The proof follows by standard arguments of optimal design theory and is therefore omitted.

\begin{satz}\label{theo:aequiv}
The design $\xi^{\star}= (\xi^{\star}_1, \ldots, \xi_M^{\star}, \mu^{\star})$ is $D$-optimal if and only if the $M$ inequalities
\be \label{eq:aequiv}
\kappa_i(d, \xi^{\star}, \theta) = h^T_i(d) M^{-1}(\xi^{\star}, \theta) h_i(d) \leq m=p+qM ,\
\ee
are satisfied for all $d \in \mathcal{X}_i$, $i= 1, \ldots, M $.
Equality holds in \eqref{eq:aequiv} for any points $(d_1, \ldots, d_M, z) \in \supp(\xi^\star_1) \times
 \ldots \times \supp(\xi^\star_M) \times \supp(\mu^\star)$.
\end{satz}

\noindent
Denote
\begin{equation} \label{xim}
\Xi^M_m = \Big\{\xi=(\xi_1, \ldots, \xi_M, \mu)~\Big| \sum_{i=1}^{M} |\supp(\xi_i)|= m \Big\}
\end{equation}
as  the set  of all designs on  $ \mathcal{X}_1 \times \ldots \times  \mathcal{X}_M \times \{1,\ldots , M \}$  with
exactly  $m$ different dose levels in the $M$ groups. The proof of the next   lemma  follows by similar arguments as
in the standard case [see \cite{silvey1980} among others], and is therefore also
omitted.

\begin{lem}\label{lem:doptweights}
Let  $\xi  =(\xi_{1}, \ldots, \xi_{M}, \mu)  \in \Xi^M_m $
denote  a  design on $ \mathcal{X}_1 \times \ldots \times  \mathcal{X}_M \times \{1,\ldots , M \}$  and $m_i$ denote the number of support points of $\xi_i$ $(i=1,\ldots, M)$. Assume that the $m= \sum^M_{i=1}m_i$ vectors \linebreak
$h_1(d^{(1)}_1), \ldots, h_1(d^{(1)}_{m_1}), \ldots, h_M(d^{(M)}_1), \ldots, h_M(d^{(M)}_{m_M})$ are linearly independent where $d^{(i)}_j \in \supp(\xi_i)$ $j=1, \ldots, m_i$, $i=1, \ldots, M$.\\
If  $\xi$  is locally $D$-optimal in the
class $\Xi^M_m$, then each component $\xi_{i}$ has equal weights at its support points.
Moreover, the weights of $\mu$  at the points $1, \ldots, M$ are given by $\tfrac{m_1}{m} , \ldots , \tfrac{m_M}{m}$,
respectively.
\end{lem}

\noindent
  In the following two sections we
present some  locally $D$-optimal designs for the Emax, the  exponential and the  linear-in-log model. The  proofs of these
results are complicated and therefore deferred to  
Section \ref{sec:appendix}.

\subsection{Models with the same location parameter}

First, we consider the case  where only the location parameter is the same in the different models. In applications this reflects  the situation of a common placebo effect   for all groups (cf. the first column of Table \ref{tab:consmod}), and we are able to identify the locally $D$-optimal design explicitly.
We begin with a general result for  the regression functions of the form \eqref{eq:loc}
where the unknown parameter vector is given by $\theta= (\theta_1, \theta^{(1), T}_2, \ldots, \theta^{(M),T}_2)^T \in \R^{m}$ with $m= 1+ Mq$.
The following result provides a solution of the $D$-optimal design problem if the $D$-optimal design for the single models are known.
 \begin{satz}\label{theo:doptplac}
Let $\sigma^2_1= \min_{i=1, \ldots, M} \sigma^2_i$ and consider  the model  given by \eqref{eq:loc}, which satisfies
\be \label{eq:placassump}
f_0(0, \theta^{(i)}_2)=0 ~,~~
 \  \tfrac{\partial}{\partial \theta^{(i)}_2} f_0 (d, \theta^{(i)}_2) |_{d=0} = 0_q
\ee
$(i=1, \ldots, M)$. If  the design
\be \label{eq:doptsingle}
\tilde \xi^{(i)} = \begin{pmatrix} 0 & d_1^{(i)} & \ldots & d_q^{(i)} \\ \tfrac{1}{q+1} &  \tfrac{1}{q+1} & \ldots &\tfrac{1}{q+1} \end{pmatrix}
\ee
is locally $D$-optimal for the single model $f(d, \theta_1, \theta^{(i)}_2)$  ($i=1, \ldots, M$), then the locally $D$-optimal design for model \eqref{eq:loc} is given by  $\xi^{\star} = (\xi^\star_1, \ldots, \xi^\star_M, \mu^\star)$ where
\be \label{eq:Doptplac}
\xi^{\star}_1= \tilde \xi^{(1)}, \quad \xi^{\star}_i = \begin{pmatrix} d^{(i)}_1& \ldots & d^{(i)}_q \\ \tfrac{1}{q}  & \ldots & \tfrac{1}{q}\end{pmatrix}, i= 2, \ldots, M , \quad \mu^{\star}= \begin{pmatrix} 1 & 2 & \ldots &M  \\ \tfrac{q+1}{m}   & \tfrac{q}{m} & \ldots &  \tfrac{q}{m}  \end{pmatrix}.
\ee
\end{satz}

\bigskip

\noindent
Using    Theorem \ref{theo:doptplac} the placebo effect $\theta_1$ is estimated in the group where the variance is smallest (see also Lemma \ref{lem:placebo} and Corollary \ref{kor:garza}).
Moreover,  it follows from the proof of Lemma \ref{lem:placebo} that the $D$-optimal design given by Theorem \ref{theo:doptplac} is not unique if there
exist two groups, say $j^*_1$ and $j_2^*$, with $\sigma_{j^*_1}^2=\sigma_{j^*_2}^2 = \min^M_{j=1} \sigma^2_j$.
We now use these results to determine $D$-optimal designs for the Emax, exponential and linear-in-log model explicitly.

\begin{kor}\label{theo:best5point}
Let $\sigma^2_1= \min_{i=1, \ldots, M} \sigma^2_i$. The locally $D$-optimal design  for the Emax, exponential and linear-in-log
model \eqref{eq:Emax2} is of the form  $\xi^\star =~(\xi^\star_1,\ldots, \xi^\star_M,\mu^\star)$, where
\beo \label{eq:best5point}
\xi^{\star}_1= \begin{pmatrix} 0 & x^{\star, (1)} & d^{(1)}_{\max} \\ \tfrac{1}{3}  & \tfrac{1}{3}   &  \tfrac{1}{3}   \end{pmatrix},  \quad \xi^{\star}_i = \begin{pmatrix}x^{\star, (i)}& d^{(i)}_{\max} \\\tfrac{1}{2}  & \tfrac{1}{2}  \end{pmatrix}, i= 2, \ldots, M , \quad \mu^{\star}= \begin{pmatrix} 1 & 2 & \ldots &M  \\ \tfrac{3}{m}   & \tfrac{2}{m} & \ldots &  \tfrac{2}{m}  \end{pmatrix}.
\eeo
and the point $x^{\star, (i)}$ is given by
\be \label{xemax}
x^{\star, (i)} = x^{\star, (i)}_{emax}= \frac{\vartheta^{(i)}_2 d^{(i)}_{\max} }{d^{(i)}_{\max} + 2\vartheta^{(i)}_2}, \quad (i=1, \ldots, M)
\ee
for the Emax model, by
 \begin{equation}\label{xexp}
x^{\star, (i)} =   x^{\star, (i)}_{\exp} = \frac{\big (d^{(i)}_{\max} - \vartheta^{(i)}_2\big)\exp \big({d^{(i)}_{\max}}/{\vartheta^{(i)}_2}\big) + \vartheta^{(i)}_2}{\exp
\big({d^{(i)}_{\max}}/{\vartheta^{(i)}_2}\big)  -1},
  ~ (i=1, \ldots, M)
 \end{equation}
 for the exponential model and by
 \begin{equation}\label{xloglin}
x^{\star, (i)} =   x^{\star, (i)}_{\log} = \frac{\bigl(d^{(i)}_{\max} + \vartheta^{(i)}_2\bigr) { \vartheta^{(i)}_2}\log\bigl({d^{(i)}_{\max}/\vartheta^{(i)}_2 } +1 \bigr) - \vartheta^{(i)}_2 d^{(i)}_{\max}}{d^{(i)}_{\max}}, ~ (i=1, \ldots, M)
 \end{equation}
 for the linear-in-log model.
\end{kor}

\noindent
It is worthwhile to mention that the locally $D$-optimal design for   model \eqref{eq:loc} with an Emax curve consists of the designs which are locally $D$-optimal for the models given by an individual Emax model with parameter $(\theta_1, \theta^{(1)}_2)$ and by an Emax model with location parameter equal to zero and parameter $\theta^{(i)}_2$, $i=2, \ldots, M$. This effect can also be observed for the exponential and the linear-in-log model.

\subsection{Models with the same location and scale parameters}

In this section we consider model \eqref{locscal} and assume that the location and scale parameter coincide across the different models
(cf. the second column in Table \ref{tab:consmod}). It turns out that in this case the $D$-optimal design problem is substantially harder, and for the
sake of a transparent presentation, we restrict ourselves to the case of $M=2$ groups. Similar results can be obtained in the case $M > 2$ with an additional amount of notation.
We begin  with some general properties of   locally $D$-optimal designs for the model \eqref{locscal} in the case of an   Emax,   linear-in-log and   exponential curve. For this purpose we define
$$
r = \frac {\sigma^2_1}{\sigma^2_2}
$$
as the ratio of the two population variances.
\begin{lem}\label{theo:suppoints} ~
\begin{itemize}
\item[(A)] The locally $D$-optimal design $\xi^\star= (\xi^\star_1, \xi^\star_2, \mu)$  for the Emax model \eqref{eq:Emax2a} and the linear-in-log \eqref{eq:log2a} have the following properties:
\begin{itemize}
\item[(A1)]  $|\supp(\xi^\star_1)| + |\supp(\xi^\star_2) | \in  \{4, 5\}.$
\item[(A2)]   If $|\supp(\xi^\star_1)| + |\supp(\xi^\star_2)| =  5$, then  $d^{(i)}_{\max} \in \supp(\xi^\star_i) , \, i=1,2$.
\item[(A3)] If  $|\supp(\xi^\star_1)| + |\supp(\xi^\star_2)| =  4$, then $d^{(1)}_{\max} \in  \supp(\xi^\star_1)$ or $d^{(2)}_{\max} \in  \supp(\xi^\star_2)$.
\end{itemize}
\item[(B)]
The locally $D$-optimal design $\xi^\star= (\xi^\star_1, \xi^\star_2, \mu)$  for the exponential model \eqref{eq:exp2a} satisfies
$$|\supp(\xi^\star_1)| + |\supp(\xi^\star_2) | \in  \{4, 5, 6 \}.$$
\end{itemize}
\end{lem}

By the  previous lemma the number of support points of the locally $D$-optimal designs is at most
$5$ for the Emax and linear-in-log model and at most $6$ for  the exponential model. On the other hand,
at least four support points are required to estimate all parameters  in both models (note that the scale and location
are assumed to be the same throughout this section). In the following discussion we determine such
``minimally'' supported  $D$-optimal  designs explicitly  for the Emax, exponential and linear-in-log model.

\subsubsection{Minimally supported designs} \label{secmin}

Recall the definition of the set $ \Xi^M_m $  in \eqref{xim}. We call a design  of the   $\xi =(\xi_1,\xi_2,\mu)$ minimally supported (for the Emax, linear-in-log and exponential
model)  if $\xi\in  \Xi^2_4 $ (note that for these models  the information matrix is of size $4 \times 4$ as the scale and location parameter coincide in both models).
It turns out the the minimally supported $D$-optimal designs for the three models under consideration have a very similar structure. On the other hand the question, if these designs
are $D$-optimal in the class of all designs does not have a simple answer and will be   discussed in the following section.

\begin{satz}\label{theo:best4point}
Let $\bar \theta^{(i)}_2 = \tfrac{ \theta^{(i)}_2}{d^{(i)}_{\max}}$, $i=1, 2$ and $0 < \bar \theta^{(1)}_2 < \bar \theta^{(2)}_2 < 1$, define
  $ y^\star = \theta^{(2)}_2 $,  $ z^\star = \theta^{(1)}_2 $
and $ x^{\star, (i)}  = x^{\star, (i)}_{emax}  $ by   \eqref{xemax} $(i=1,2)$.
\begin{itemize}
\item[(1)]  If $r \leq 1$,  the locally $D$-optimal design for model \eqref{eq:Emax2a}  in the class  $\Xi^2_4$
is given by
\be \label{eq:best4pointa}
\xi^{a, \star}_1= \begin{pmatrix} 0 &x^{\star, (1)}& d^{(1)}_{\max} \\  \tfrac{1}{3}  &  \tfrac{1}{3}&   \tfrac{1}{3} \end{pmatrix},  \quad \xi^{a, \star}_2 = \begin{pmatrix} y^\star
 \\ 1\end{pmatrix}, \quad \mu^{a, \star} = \begin{pmatrix} 1 & 2 \\  \tfrac{3}{4} &  \tfrac{1}{4} \end{pmatrix}.
\ee
\item[(2)]  If $ 1< r  \leq \Bigl( \tfrac{1+ \bar \theta^{(2)}_2}{1+ \bar \theta^{(1)}_2} \Bigr)^6$,  the locally $D$-optimal design for model \eqref{eq:Emax2a}  in the class $\Xi^2_4$
 is given by
\be  \label{eq:best4pointb}
\xi^{b, \star}_1= \begin{pmatrix} x^{\star, (1)} & d^{(1)}_{\max} \\ \tfrac{1}{2}  & \tfrac{1}{2}  \end{pmatrix},  \quad \xi^{b, \star}_2 = \begin{pmatrix} 0 &  y^\star  \\ \tfrac{1}{2} & \tfrac{1}{2} \end{pmatrix}, \quad \mu^{b, \star} = \begin{pmatrix} 1 & 2 \\ \tfrac{1}{2}  & \tfrac{1}{2}  \end{pmatrix}.
\ee
\item[(3)]   If $r >  \Bigl( \tfrac{1+ \bar \theta^{(2)}_2}{1+ \bar \theta^{(1)}_2} \Bigr)^6$, the locally $D$-optimal design for model \eqref{eq:Emax2a}  in the class  $\Xi^2_4$
 is given by
\be  \label{eq:best4pointc}
\xi^{c, \star}_1 = \begin{pmatrix}z^\star \\ 1\end{pmatrix}, \quad \xi^{c, \star}_2= \begin{pmatrix} 0 & x^{\star, (2)} & d^{(2)}_{\max} \\  \tfrac{1}{3}  &  \tfrac{1}{3}&   \tfrac{1}{3}\end{pmatrix},  \quad \mu^{c, \star} = \begin{pmatrix} 1 & 2 \\  \tfrac{1}{4} &  \tfrac{3}{4} \end{pmatrix}.
\ee
\end{itemize}
\end{satz}

We can also obtain the  minimally supported $D$-optimal designs for the exponential model and the linear-in-log with common location and scale
parameter.

\begin{satz}\label{theo:best4pointexp}
Let $\bar \theta^{(i)}_2 = \tfrac{ \theta^{(i)}_2}{d^{(i)}_{\max}}$, $i=1, 2$,  $0 < \bar \theta^{(1)}_2 < \bar \theta^{(2)}_2 < 1$,
define
 $$g(\theta, x) = \left( 1 +(x-1)  \exp(\tfrac{x}{\theta})- x\exp(\tfrac{x-1}{\theta}) \right)^2$$
and   $y^{\star}=  d^{(2)}_{\max} $,  $z^{\star}=  d^{(1)}_{\max} $ and   the point $x^{\star,(i)}= x_{\exp}^{\star,(i)}$ by \eqref{xexp} for $i=1, 2$.

\begin{itemize}
\item[(1)]  If $r \leq 1$,  the $D$-optimal design for model \eqref{eq:exp2a}  in the class  $\Xi^2_4$
is given by \eqref{eq:best4pointa}.
\item[(2)]  If $ 1< r  \leq  \tfrac{g(\theta^{(1)}_2, x^{\star, (1)}_{\exp})}{g(\theta^{(2)}_2,  x^{\star, (2)}_{\exp})}$,  the $D$-optimal design for model \eqref{eq:exp2a}  in the class $\Xi^2_4$
 is given by \eqref{eq:best4pointb}.
\item[(3)]  If $r >  \tfrac{g(\theta^{(1)}_2,  x^{\star, (1)}_{\exp})}{g(\theta^{(2)}_2,  x^{\star, (2)}_{\exp})}$, the $D$-optimal design for model \eqref{eq:exp2a}  in the class  $\Xi^2_4$
 is given by \eqref{eq:best4pointc}.
\end{itemize}
\end{satz}

\medskip

\begin{satz}\label{theo:best4pointlog}
Let $\bar \theta^{(i)}_2 = \tfrac{ \theta^{(i)}_2}{d^{(i)}_{\max}}$ $i=1, 2$,  $0 < \bar \theta^{(1)}_2 < \bar \theta^{(2)}_2 < 1$, define 
 $$g(\theta, x) = (1 + \theta )^2 \big (\log(\tfrac{1}{\theta} +1 )  \log(\tfrac{x}{\theta}+1 ) \big)^2 \Big( \frac{x}{ (x+\theta) \log(\tfrac{x}{\theta} +1 ) } - \frac{1}{(1+\theta) \log(\tfrac{1}{\theta} +1 ) } \Big )^2$$
and $y^{\star}=  d^{(2)}_{\max} $, $z^{\star}=  d^{(1)}_{\max} $
 and   the point $x^{\star,(i)}= x_{\log}^{\star,(i)}$ by  \eqref{xloglin} for $i=1, 2$.
\begin{itemize}
\item[(1)]  If $r \leq 1$,  the $D$-optimal design for model \eqref{eq:log2a}  in the class  $\Xi^2_4$
is given by \eqref{eq:best4pointa}.
\item[(2)]  If $ 1< r  \leq  \tfrac{g(\theta^{(1)}_2, x_{\log}^{\star,(1)})}{g(\theta^{(2)}_2, x_{\log}^{\star,(2)})}$,  the $D$-optimal design for model \eqref{eq:log2a}  in the class $\Xi^2_4$
 is given by \eqref{eq:best4pointb}.
\item[(3)]  If $r >   \tfrac{g(\theta^{(1)}_2, x_{\log}^{\star,(1)})}{g(\theta^{(2)}_2, x_{\log}^{\star,(2)})} $, the $D$-optimal design for model \eqref{eq:log2a}  in the class  $\Xi^2_4$
 is given by \eqref{eq:best4pointc}.
\end{itemize}
\end{satz}

\subsubsection{$D$-optimal designs in the class of all designs}

The question if a minimally supported $D$-optimal design for one of the models considered in Section \ref{secmin}
is in fact $D$-optimal  in the class of all designs is an extremely difficult one. Its answer depends sensitively on the
particular parameters in the model under  consideration  and  differs for the three dose response models under consideration.
We exemplarily state a result for the Emax model, which  provides sufficient conditions for the $D$-optimality of a minimally supported
$D$-optimal design, and illustrates the general
structure and difficulties in results of this type. The proof is based on the equivalence Theorem \ref{eq:aequiv} and given in the appendix. Similar but substantially
more complicated statements  can  also be obtained of the linear-in-log and the exponential model (note that in contrast to the Emax model
these models contain transcendental functions).

\begin{satz}\label{theo:Emaxdopt}
 Let $\bar \theta^{(i)}_2 = \tfrac{ \theta^{(i)}_2}{d^{(i)}_{\max}}$, $i=1, 2$ and assume $0 < \bar \theta^{(1)}_2 < \bar \theta^{(2)}_2 < 1$.
 \begin{itemize}
 \item[(1)]  Let $r \leq 1$. The design $\xi^{a, \star}$ defined in
  \eqref{eq:best4pointa} is locally $D$-optimal for model \eqref{eq:Emax2a} if the condition
\be \label{eq:Emaxdopt4a1}  
\bar \theta^{(2)}_2 \geq \frac{ r \bigl( 6 \bar \theta^{(1)}_2 (\bar \theta^{(1)}_2+ 1) (2\bar \theta^{(1)}_2 +1)^2 \bigr) -
 \bigl(1-  r \bigr) }{(6 + 2 r \bar \theta^{(1)}_2(1+ 2\bar \theta^{(1)}_2))}
\ee
is satisfied.
 \item[(2)]  Let  $r > 1$. The design $\xi^{b, \star}$
defined in  \eqref{eq:best4pointb} is locally $D$-optimal for model \eqref{eq:Emax2a} if and only if the condition
\be \label{eq:Emaxdopt4b} 
 \bar \theta^{(2)}_2 \geq \frac{   (\bar \theta^{(1)}_2)^2 (1+ 2 \bar \theta^{(1)}_2)^2 + r (1+ \bar \theta^{(1)}_2)^2 (1+ 4 \bar \theta^{(1)}_2 + 20 (\bar \theta^{(1)}_2)^2) -1}{6 + 2 \bar \theta^{(1)}_2(1+ 2\bar \theta^{(1)}_2)}
\ee
is satisfied.
 \item[(3)] Let $r>1$. The design $\xi^{c, \star}$
defined in \eqref{eq:best4pointc} is locally $D$-optimal for model \eqref{eq:Emax2a} if the condition
\be \label{eq:Emaxdopt4c1}  
\bar \theta^{(1)}_2 \geq \frac{ \frac{1}{r} \Bigl( 6 \bar \theta^{(2)}_2 (\bar \theta^{(2)}_2+ 1) (2\bar \theta^{(2)}_2 +1)^2 \Bigr) - \Bigl(1-   \frac{1}{r} \Bigr) }{(6 + 2 \frac{1}{r} \bar \theta^{(2)}_2(1+ 2\bar \theta^{(2)}_2))} 
\ee
is satisfied.
\end{itemize}
\end{satz}

\begin{figure}[t!]
\centering
 \includegraphics[width=0.95\textwidth]{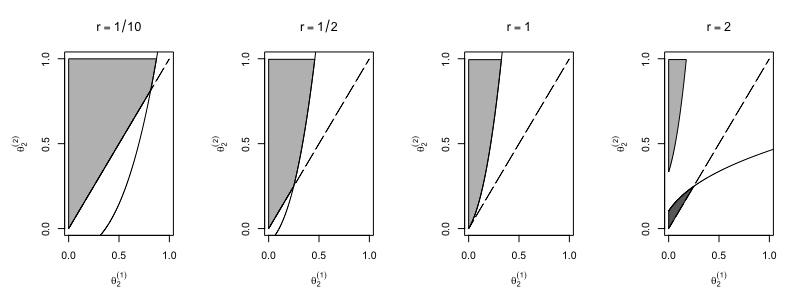}
\caption{\label{Fig:domain}
\it The marked regions describe the parameter spaces, where the minimally supported $D$-optimal design 
is  optimal in the Emax model  \eqref{eq:Emax2a} (see Theorem  \ref{theo:Emaxdopt}). The different figures correspond to different values of $r={\sigma^2_1}/{\sigma^2_2}$. 
The domain for the first case of  Theorem~\ref{theo:Emaxdopt} is represented in gray for 
the case $r=1/10$, $r=1/2$ and $r=1$ (see the first three  panels from the left). In the right  panel we display the case $r=2$  of  Theorem~\ref{theo:Emaxdopt}
(here the  gray region  corresponds to case (2), while the dark gray region corresponds to case (3)).}
\end{figure}

Figure \ref{Fig:domain} illustrates the parameter domains for different ratios $r= \tfrac{\sigma^2_1}{\sigma^2_2}$. The case where the variance is equal in both groups is presented in the third panel. Obviously, there are several  parameter constellations $\theta^{(1)}_2 \geq \theta^{(2)}_2$ where the minimally supported $D$-optimal design $\xi^{a, \star}$ is not optimal in the class of all designs. \\


\section{Application to a dose-finding study}\label{secEX2}
\def\theequation{5.\arabic{equation}}
\setcounter{equation}{0}

In this section we illustrate the application of the results of the previous sections
and discuss the problem of designing experiments for a dose finding study with different treatment groups.
Our example refers to a Phase II study on a drug that works by
increasing the level of a biomarker that induces a beneficial clinical effect in patients.  The dosing
groups under consideration are monthly and weekly
administration. The primary objective of the study is the
characterization of dose-response relationships at a given
time-point, say  $T$, after initiation of treatment for each of these two dosing
groups.  This will support the selection of an appropriate dose level
and group to be used in phase III clinical trials. To maintain the confidentiality
of the trial the dose-range has been rescaled and the considered range
(in terms of total monthly dose) is $[0,400]$ for the weekly group
and $[0,1000]$ for the once-a-month group.
The natural questions for the  design of this study are (i)
 which doses should be studied  in each treatment
group and (ii) how to split the total sample size between the two treatment groups.
Here the objectives of the
  study are addressed by deriving the best estimates
  of the dose response curves, a task for which a $D$-optimal design
  is best suited.

To arrive at a suitable design for the Phase II study, we need to
quantify the information. This quantification can generate
 a best guess for the dose-response curves, but, even better, it can
 be used to obtain a candidate set of dose-group-response scenarios
  to reflect the uncertainty about the true
dose-group-response relationship.
The available information was data from a very small early trial,
which was used to develop a nonlinear mixed
effects pharmacokinetics (PK) / pharmacodynamics (PD) model linking
drug concentrations to biomarker levels. Using this model, data of the new trial were predicted for
the time-point $T$ of the dose-response analysis and
dose-group-response models were fitted to the data. Under the
assumption of a normal distribution for the logarithm of the biomarker level, 
it turned out that the Emax function was able to adequately describe
the population average predicted by the PK/PD model. The Emax model
utilized total monthly dose as input and had different $\mbox{ED}_{50}$
parameters in the two groups ($\theta_2^{(1)}$ and $\theta_2^{(2)}$),
but the same placebo $\vartheta_{11}$ and Emax parameter
$\vartheta_{12}$, so that the model function in the weekly and monthly
group is given by
\beo
f (d, \theta_1, \theta^{(i)}_2)= \vartheta_{11} + \vartheta_{12}
\frac {d}{\theta^{(i)}_2+d}, \quad i=1, 2,
\eeo
Here group $i=1$ contains patients receiving monthly administration
and the group $i=2$ the weekly administration. The parameter
estimates can be found in Table \ref{tbl:cand-set}  as model 1, which can be considered as population average fit.
We now use these estimates as a guess 
and determine the locally $D$-optimal design for these values.   The variability is
expected to be the same in both
treatment groups. Recalling the  design spaces for
the monthly and weekly doses are $\mathcal{X}_1=[0,1000]$ and
$\mathcal{X}_2=[0,400]$, respectively, we obtain from 
 Theorem~\ref{theo:best4point} and
Theorem~\ref{theo:Emaxdopt} the (locally) 
$D$-optimal design $\xi^\star=(\xi_{1}^\star,\xi_{2}^\star,\mu^\star)$
 as
\begin{eqnarray*}
\xi_{1}^\star= \begin{pmatrix} 0 &x^{\star, (1)}_{emax} & d^{(1)}_{\max} \\  \tfrac{1}{3}  &  \tfrac{1}{3}&   \tfrac{1}{3} \end{pmatrix} = \begin{pmatrix} 0 &13.45  &	1000 \\  \tfrac{1}{3}  &  \tfrac{1}{3}&   \tfrac{1}{3} \end{pmatrix} ,
\quad \xi_{2}^\star = \begin{pmatrix} \theta^{(2)}_2 \\ 1\end{pmatrix} = \begin{pmatrix} 10.46\\ 1\end{pmatrix} ,
\mu^\star = \begin{pmatrix} 1 & 2 \\  \tfrac{3}{4} &  \tfrac{1}{4} \end{pmatrix}.
\end{eqnarray*}
It can be seen that based on the population average fit, it is
sufficient to investigate the low dose-range in both groups and a
high dose in one of the two groups. Here the maximum dose is placed in
the monthly group because
$\theta^{(1)}_2 /d^{(1)}_{\max}< \theta^{(2)}_2 / d^{(2)}_{\max}$, so
relative to the allowed maximum dose a larger $\mbox{ED}_{50}$ parameter
exists for the weekly group and thus patients are allocated to the
monthly group.

\begin{figure}[t]
\centering
 \includegraphics[width=0.65\textwidth]{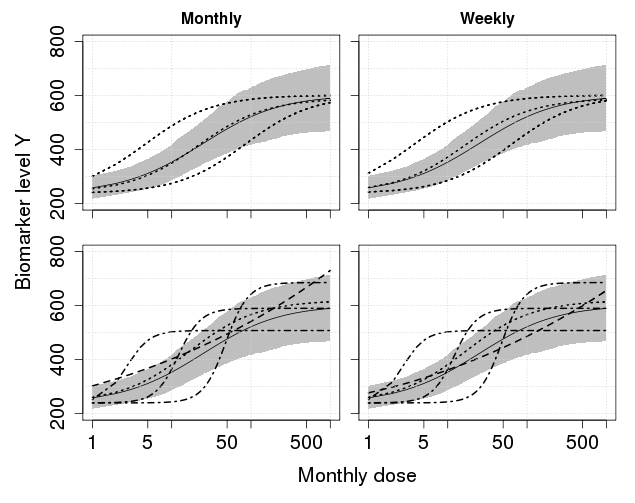}
\caption{\it {Candidate models for the dose-response curve in monthly and
  weekly group. In the first row models 1-5 are depicted and in the
  second row models 6-10 (see Table \ref{tbl:cand-set}). The solid line represents the population
  average of the new trial data generated with the PK/PD model. The grey area represents the biomarker level between the 25th and 75th quantiles of the patient responses  of the new trial data. Dotted curves correspond to the Emax models.
  Sigmoid Emax models are depicted as dotted-dashed lines
  and the linear-in-log models as dashed lines.}}
\label{Fig:models}
\end{figure}

 In practice, it is not realistic to assume that the data and model 
from previous small trials completely represent the underlying truth
(otherwise no further study would need to be conducted). So it is important to derive ranges covering the uncertainty
about the available information to use for the design of the new
study. In particular, because the population to be included in the
Phase II trial will cover a broader range of characteristics than in the small proof
of concept trial.  For this purpose the PK/PD model was used to
predict individual dose-response curves and the Emax model was  
fitted to the individual dose-response curves to derive a range of
plausible dose-response parameters. Quantiles of the derived parameter
distributions were used to derive four additional candidate model
shapes. More details on how these candidate shapes were derived
can be found in Appendix \ref{sec:appendix2}. The parameters for these
four additional candidate models can be found in Table
\ref{tbl:cand-set} under the numbers 2-5. These models are depicted
in the first row of  Figure \ref{Fig:models}.

\begin{table}[t]
\begin{center}
\begin{tabular}{|c|c|cccccc|}
\hline
model type & id & $\vartheta_{11} $ & $\vartheta_{12} $ & & $\theta_2^{(1)}$ & $\theta_2^{(2)}$ & $\gamma$  \\
\hline
Emax & 1 & 5.48& 0.90& & 13.82& 10.46 &    1  \\
Emax &2 & 5.47 &  0.93 &   &  2.93 &  2.39 &   1 \\
Emax &3 &5.47 &  0.93 &  &  2.93 & 40.40 &   1 \\
Emax &4 & 5.47 & 0.93 &  & 53.49 &  2.39 &   1 \\
Emax &5 & 5.47 &  0.93 &  & 53.49 & 40.40 &    1 \\
Sigmoid Emax &6 & 5.48 & 0.90 &  & 13.82 & 10.46 &    3 \\
\hline
Model type & id  & $\theta_1 $ & $\vartheta_1^{(1)}$ & $\vartheta_1^{(2)}$     & $\vartheta_2^{(1)}$ & $\vartheta_2^{(2)}$ & $\gamma$    \\
\hline
Emax &7 & 5.48 &  0.85 &  0.95 &  13.82 & 10.46 &    1 \\
Sigmoid Emax &8 & 5.48 & 0.65 & 0.75 &  2.93 &  2.39 &    3 \\
Sigmoid Emax &9 & 5.48 & 0.95 & 1.05 & 53.49 & 40.40  &    3 \\
Log & 10 & 5.44 &  0.13 & 0.14 &  0.32 &  0.41 &     \\
\hline
\end{tabular}
\caption{\it Set of candidate models used in the robust criterion \eqref{compound}}
\label{tbl:cand-set}
\end{center}
\end{table}

 With this set of candidate
models, the design maximizing the mean 
efficiency 
\be \label{compound}
g_c(\xi,s)=\sum_{i=1}^{s} \pi_i \mbox{Eff}_i(\xi)~
\ee
can be calculated,
where $s$ is the number of candidate models (here $5$ or $10$),
$\pi_1,\ldots,\pi_{s}$ are nonnegative model weights chosen to reflect
prior probability associated the model function $1,\ldots,s$ (throughout this paper we
will use $\pi_i=1/s$, $i=1,\ldots,s$). The efficiencies
$\mbox{Eff}_i(\xi)$ of the experimental design $\xi$ with respect to the (locally)
$D$-optimal design $\xi^{\star,i}$ associated to the model $i$ is defined as
$$\mbox{Eff}_i(\xi )=\left(\frac{|M_i(\xi,\theta_i))|}{|M_i(\xi^{\star,i},\theta_i)|}\right)^{1/m_i},$$
where $M_i$ is the Fisher information matrix associated to the model $i$ with parameter specification $\theta_{i}$
and $m_i$ is the number of parameters of this model. The criterion \eqref{compound} is called {\it Bayesian} 
or {\it  compound} optimality criterion in the literature [see \cite{dette1990}, \cite{cookwong1994}
or \cite{tsaizen2004,zentsa2004} among many others].
In the following we will
denote the designs maximizing  (\ref{compound}) by $\xi^\star_{c,s} = (\xi^\star_{1,c,s},\xi^\star_{2,c,s},\mu_{c,s}^\star)$
and call it compound optimal design. We emphasize that the definition of the criterion \eqref{compound} 
requires knowledge  of the locally optimal designs $\xi^{\star,i}$, which have been determined in Section \ref{secRES}.
\\
The  compound optimal design based on the first $5$ models  in Table
\ref{tbl:cand-set} can be  calculated numerically and is given by
 $\xi^\star_{c,5}=(\xi^\star_{1,c,5},\xi^\star_{2,c,5},\mu_{c,5}^\star),$ where
\begin{eqnarray*}
\xi_{1, c,5}^\star \approx \begin{pmatrix} 0		& 3.02 	&   43.67	& 1000 \\  0.26  & 0.24&0.25 &0.25   \end{pmatrix} ,
\quad \xi_{2, c,5}^\star = \begin{pmatrix}  2.53	& 37.51\\ 0.48 & 0.52\end{pmatrix} ,
\mu_{c,5}^\star = \begin{pmatrix} 1 & 2 \\ 0.67 &  0.33\end{pmatrix}~,
\end{eqnarray*}
and its optimality can be proved by an analogue of Theorem \ref{theo:aequiv} for the Bayesian optimality criterion \eqref{compound}.
Compared to the design using only the best guess model, now the low dose-range is investigated in finer
granularity by using two instead of one dose (safeguarding
against different possible values of the $\mbox{ED}_{50}$).
In addition still more patients are evaluated for the monthly
group, as the high dose is only used there.

 Based on general
plausibility considerations five further candidate shapes were included
as example of models different from the Emax function (\textit{e.g.} the
sigmoid Emax and the linear-in-log function), or of models where the
maximum efficacy differed between the two groups. These   models
are shown  in the second row  of  Figure \ref{Fig:models} and the corresponding
parameters are given  in the  rows  with ids 6-10 in Table \ref{tbl:cand-set}.
 First a sigmoid Emax model
 \beo \label{eq:SigEmax1}
 f (d, \theta_1, \theta^{(i)}_2,\gamma)= \vartheta_{11} + \vartheta_{12} \frac {d^{\gamma}}{(\theta^{(i)}_2)^{\gamma}+d^{\gamma}}, \quad i=1, 2,
\eeo 
with Hill coefficient $\gamma=3$  (model 6) is also considered as a possible
dose response function. Note that  this  model provides a steeper
dose-response curve compared to the Emax model, but with the same
$\mbox{ED}_{50}$ values as model 1. Furthermore, an Emax and a sigmoid Emax model 
\beo \label{eq:SigEmax2}
 f (d, \theta_1, \theta^{(i)}_2,\gamma)= \theta_{1} + \frac {\vartheta^{(i)}_{1}  d^\gamma}{(\vartheta^{(i)}_2)^{\gamma}+d^{\gamma}}, \quad i=1, 2.
\eeo
is added that allows for
different Emax parameters in the two treatment group (models 7, 8, 9).
In addition a linear-in-log model (id 10) is utilized. The locally  $D$-optimal designs
for these models  can be computed  using the results of Section~\ref{secRES}.
For  the sigmoid Emax models, a transformation has  to be used 
 to reduce it to the case of an Emax model, such that the derived theory is applicable
 (note that the parameter  $\gamma$ is assumed to be fixed).

When using all $s=10$ candidate models we obtain $\xi^\star_{c,10}=(\xi^\star_{1,c,10},\xi^\star_{2,c,10},\mu_{c,10}^\star)$ where
\begin{eqnarray*}
\xi_{1, c,10}^\star \approx \begin{pmatrix} 0		& 2.90	& 12.98   &41.91	& 1000  \\ 
 0.27  & 0.13 &0.22 &0.13 & 0.24   \end{pmatrix} ,
\quad \xi_{2, c,10}^\star = \begin{pmatrix}  3.01	& 13.16	&49.46	& 400\\  0.33&0.21&0.31 & 0.15\end{pmatrix} ,
\end{eqnarray*}
$$\mu_{c,10}^\star = \begin{pmatrix} 1 & 2 \\ 0.58 &  0.42\end{pmatrix}.$$
This design investigates the lower dose range comparably to the
previous design based on the first five candidate models, but the
maximum dose is studied in both groups.
\begin{table}
\begin{tabular}{|c|c|cccccccccc|}
\hline
&$g_c(\cdot, s)$ & 1 & 2 & 3 & 4 & 5 & 6 & 7 & 8 & 9 & 10  \\
\hline
$\xi^\star_{c,5}$  &  0.823 &	0.708 &0.835 &0.877 &0.845 &0.847 & 0.098 &0.795 &0.927 &0.906 &0.625	\\
$\xi^\star_{c,10}$  &0.747 & 0.831 &0.749 &0.779 &0.767 &0.786 &0.749 &0.903 &0.760 &0.749 & 0.747	\\
\hline
\end{tabular}
\caption{\it Efficiency $\mbox{Eff}_i(\xi^\star_{c,s} )$ of the two compound optimal designs compared to each of the locally $D$-optimal designs for the 10 models.}
\label{tbl:eff}
\end{table}
The efficiencies of the two  designs $\xi^\star_{c,5}$ and $\xi^\star_{c,10}$  in the different models are displayed in Table \ref{tbl:eff}.
We observe that the design $\xi^\star_{c,5}$  has reasonable efficiencies in all models except in the  sigmoid Emax (6). Note that this
design has been constructed on the basis of the models (1) - (5). On the other hand the
 the design $\xi^\star_{c,10}$ maximizes the criterion \eqref{compound}, where  uncertainty with respect to all models (1) - (10) is addressed.
 As a consequence this design has efficiencies varying between $75\%$ - $90\%$ in all  competing models under consideration. 
Moreover,  it can be used for a goodness-of-fit test of the Emax model, as both components have more than $3$ support points.
For these reasons we recommend this design for the Phase II study under consideration.

\bigskip

{\bf Acknowledgements} The authors would like to thank  Antoine Soubret for the PK/PD model
and Martina
Stein, who typed parts of this manuscript with considerable
technical expertise. We are also grateful to Katrin Kettelhake for computational assistance
and to Antoine Soubret, who built the original PK/PD model that was used
to derive the candidate models for the designs calculated in Section \ref{secEX2}.
This work has been supported in part by the
Collaborative Research Center "Statistical modeling of nonlinear
dynamic processes" (SFB 823, Teilprojekt C2) of the German Research Foundation
(DFG) and by a grant from the National Institute Of General Medical Sciences of the National
Institutes of Health under Award Number R01GM107639. 
The content is solely the responsibility of the authors and does not necessarily
 represent the official views of the National
Institutes of Health.

\begin{small}
\bibliographystyle{apalike}
\bibliography{chrystel_lit}
\end{small}

\bigskip
\medskip

\section{Appendix: technical details} \label{sec:appendix}
\def\theequation{6.\arabic{equation}}
\setcounter{equation}{0}

\subsection{Deviation of the information matrix}  \label{A1}
Assuming a normal distribution of the errors  $(\varepsilon_{ij\ell} \sim \mathcal{N}(0,\sigma^2_i)$ independent) 
partial derivatives of the corresponding log-likelihood function $\ell(\theta)$  with respect to $ \theta_1$ and $\theta^{(i)}_2$ are given by
 \begin{eqnarray*}
 \frac {\partial}{\partial \theta_1} \ell(\theta) &=&  \sum^M_{i=1} \sum^{k_i}_{j=1} \sum^{n_{ij}}_{\ell =1} \frac {1}{\sigma^2_i}  \Bigl(Y_{ij \ell}   -  f(d^{(i)}_j, \theta_1, \theta^{(i)}_2))\Bigr)  \frac{\partial}{\partial \theta_1} f(d^{(i)}_j, \theta_1, \theta^{(i)}_2), \\
 \frac {\partial}{\partial \theta^{(i)}_2} \ell(\theta) &=&   \sum^{k_i}_{j=1} \sum^{n_{ij}}_{\ell =1} \frac {1}{\sigma^2_i}  \Bigl(Y_{ij \ell}   -f(d^{(i)}_j, \theta_1, \theta^{(i)}_2) \Bigr) \frac {\partial}{\partial \theta^{(i)}_2} f(d^{(i)}_j, \theta_1, \theta^{(i)}_2),
 \end{eqnarray*}
 ($i= 1, \ldots, M$).  Note that
 \begin{eqnarray*}
 \mathbb{E} \Bigl[ \frac {\partial}{\partial \theta^{(i)}_2} \ell (\theta) \Bigl( \frac {\partial}{\partial \theta^{(i)}_2} \ell (\theta) \Bigr)^T \Bigr] &=&
 \left \{
 \begin{array} {cc}
 0 & \mbox{if} \ \ i \neq i^\prime \\
 \sum^{k_i}_{j=1} \sum^{n_{ij}}_{\ell=1} \frac {1}{\sigma^2_i} \eta_i (d^{(i)}_j) \eta^T_i (d^{(i)}_j) & \mbox{if} \ \ i = i^\prime
\end{array}
 \right . \\
\mathbb{E} \Bigl [ \frac {\partial}{\partial \theta_1} \ell(\theta) \Bigl( \frac {\partial}{\partial \theta^{(i)}_2} \ell(\theta) \Bigr)^T \Bigl] &=& \frac {1}{\sigma^2_i} \sum^{k_i}_{j=1} \sum^{n_{ij}}_{\ell=1} \overline \eta_i (d^{(i)}_j) \eta_i (d^{(i)}_j)
 \end{eqnarray*}
 and
 $$
 \mathbb{E} \Bigl[ \frac {\partial}{\partial \theta_1} \ell (\theta) \Bigl( \frac {\partial}{\partial \theta_1} \ell (\theta) \Bigr)^T \Bigr] = \sum^M_{i=1} \sum^{k_i}_{j=1} \sum^{n_{ij}}_{\ell=1} \frac {1}{\sigma^2_i}  \overline \eta_i (d^{(i)}_j) \overline \eta^T_i (d^{(i)}_j), \qquad \qquad
 $$
 where
 $\overline \eta_i^T (d)  =   \tfrac{\partial}{\partial \theta_1} f(d^{(i)}_j, \theta_1,  \theta^{(i)}_2),  \
  \eta_i^T (d) =  \tfrac{\partial}{\partial \theta^{(i)}_2} f(d^{(i)}_j, \theta_1, \theta^{(i)}_2)$.  
  Consequently, we obtain for the Fisher information
 $$
M_n:= \sum^M_{i=1} \sum^{k_i}_{j=1} n_{ij} h_i (d^{(i)}_j) h^T_i (d^{(i)}_j) =
 n \sum^M_{i=1} \frac {n_i}{n} \sum^{k_i}_{j=1} \frac{ n_{ij}}{n_i} h_i (d^{(i)}_j) h^T_i (d^{(i)}_j) ,
 $$
 where the vector $h_i$ is defined in \eqref{gradient}. Observing the assumption \eqref{limit} it follows that $\frac {1}{n}M_n$ converges to the matrix $M(\xi,\theta)$ defined in \eqref{eq:info}.

\subsection{Proof of main results}
{\bf Proof of Lemma \ref{theo:garza}:}
We only discuss the first part of the proof. The second assertion follows by similar argument.
Let $\xi= (\xi_1, \ldots, \xi_M, \mu)$ be an arbitrary design with $|\supp(\xi_i)| \geq \tfrac{(k+2)}{2}$ $(i=1,\ldots,M)$ and assume that there exists a design $\xi^{+} =  (\xi^{+}_1, \ldots, \xi^{+}_M, \mu)$ such that
$$C(\xi^{+}, \theta_1, \theta^{(i)}_2) \geq C(\xi , \theta_1, \theta^{(i)}_2)$$
for all $i=1, \ldots, M$. Recalling the definition of the matrix $P_i$ in \eqref{pi} it then   follows by Theorem 14.2.9 of \cite{harv1997} that
$$M^{(i)} (\xi^+, \theta) = P_i C(\xi^{+}, \theta_1, \theta^{(i)}_2) P^T_i \geq  P_i C(\xi,  \theta_1, \theta^{(i)}_2) P^T_i = M^{(i)}(\xi, \theta) \quad i=1, \ldots, M.$$
This implies
$$M(\xi^+, \theta) = \sum_{i=1}^{M}\lambda_i M^{(i)}(\xi^{+}_i, \theta) \geq  \sum_{i=1}^{M}\lambda_i M^{(i)}(\xi_i, \theta) =M(\xi, \theta)$$
and the design $\xi^+$ increases the information matrix $M( \cdot, \theta)$ with respect to the Loewner ordering. \\
It now follows from
  Theorem 3.1 of \cite{detmel2011}   that there exists a design $\xi^+$ with components $\xi^+_i$ with at most $\tfrac{k+2}{2}$ support points   ($i=1, \ldots, M$). The statements (1a) and (1b) in Lemma \ref{theo:garza} also follows from Theorem 3.1 in \cite{detmel2011}.

\bigskip

\noindent
{\bf Proof of  Lemma \ref{lem:placebo}:}
We only prove the Lemma for the model given by \eqref{locscal}. The proof for model \eqref{eq:loc} is analogous.
Note that in the  model under consideration we have
  $ \tfrac{\partial}{\partial \theta_1} f(d, \theta_1, \theta^{(i)}_2)=  \big(1,  f_0(d, \theta^{(i)}_2) \big)$ for the gradient in \eqref{gradient}. Consequently, if $\delta_0$ denotes the Dirac measure at the point $0$, it follows
for the matrices $M^{(i)}$  defined in \eqref{eq:infoi} that
\begin{equation}
\label{matplac}
\sigma_i^2 M^{(i)}(\delta_0, \theta) = \sigma_1^2 M^{(1)}(\delta_0, \theta)~, ~i= 1, \ldots, M.
\end{equation}
\noindent
Now, we consider the design $\eta = (\eta_1, \ldots, \eta_M, \nu )$ and represent its components as
\begin{equation*}
\label{eta}
\eta_i= \omega^{(i)}_0 \delta_0 + (1- \omega^{(i)}_0) \eta^0_i \, \quad i=1, \ldots, M , \qquad
\nu= \sum_{i=1}^M \lambda_i \delta_i~.
\end{equation*}
Here  $\delta_t$ is the Dirac measure  at the point $t$,     $\lambda_i, \omega^{(i)}_0 \in [0, 1]$, $i=1, \ldots, M$
 and $\eta^0_1, \ldots ,\eta^0_M$ denote designs with $0 \notin \mbox{supp}(\eta^0_i)$.
 Moreover, at least for one $i \in \{1, \ldots, M\}$  we have  $\lambda_i \omega^{(i)}_0>0$. \\
We now assume without loss of generality that $j^*=1$ and
 construct a ``better'' design  $\xi= (\xi_1, \ldots, \xi_M, \mu)$ as follows
\beo 
\label{xi}
\xi_1= \omega^{\star} \delta_0 + (1-\omega^{\star}) \eta^0_1,
\xi_i =  \eta^0_i, ~(   i=2, \ldots, M), ~
\mu= \sum_{i=1}^M \lambda^{\star}_i \delta_i  ,
\eeo
 where
\beo
 \label{eq:lambdastar}
\lambda^{\star}_1 =  \lambda_1 +  \sum_{i=2}^M \lambda_i \omega^{(i)}_0 \in [0,1] , ~
\lambda^{\star}_i  = \lambda_i (1- \omega^{(i)}_0 )~(i = 2, \ldots, M), ~
\omega^{\star} =\tfrac{\sum_{i=1}^M \lambda_i \omega^{(i)}_0 }{\lambda_1 +  \sum_{i=2}^M \lambda_i \omega^{(i)}_0}.
\eeo
Note that we shift the weights of the  measures $\eta_i$ at the point $0$ to the design for the group with the smallest population
variance.  Observing \eqref{matplac}  gives for the difference
\begin{eqnarray*}
M(\xi, \theta) - M(\eta, \theta) &=&{\omega^{\star} \lambda^\star_1}  M^{(1)}(\delta_0, \theta) + {(1-\omega^\star) \lambda^\star_1}M^{(1)}(\eta^0_1, \theta) + \sum_{i=2}^M {\lambda^\star_i}M^{(i)}(\eta^0_i, \theta) \\
&& ~- \Big(
\sum_{i=1}^{M}  {\omega^{(0)}_i \lambda_i} M^{(i)}(\delta_0, \theta) +\sum_{i=1}^M  {\lambda_i (1-\omega^{(0)}_i)} M^{(i)}(\eta^0_i, \theta)
 \Big)  \\
&=& \Big( {\omega^{\star} \lambda^\star_1}- \sum_{i=1}^M \tfrac{ {\sigma^2_1} }{\sigma^2_i} \omega^{(0)}_i \lambda_i \Big) M^{(1)}(\delta_0, \theta) \geq
\Bigl(\omega^{\star} \lambda^\star_1 - \sum_{i=1}^M \lambda_i \omega^{(i)}_0 \Bigr) \ M^{(1)}(\delta_0, \theta) = 0,
\end{eqnarray*}
since $\sigma_1^2 \leq \sigma_i^2$ ($i=1,\ldots, M$).

\bigskip
\noindent
{\bf Proof of Corollary \ref{kor:garza}:}
Lemma \ref{theo:garza} can be applied  in the case of an Emax model or linear-in-log model with $k=4$ [see \cite{yang2010}]. Consequently, there exists a design $\xi^+$ with $3M$ support points and each component $\xi^+_i$ contains the placebo $0$ and $d^{(i)}_{\max}$ $i=1, \ldots, M$. Now we apply Lemma \ref{lem:placebo} with $\eta= \xi^+$ and we allocate the placebo $0$ in the group with the smallest variance.
For the exponential model Lemma \ref{theo:garza} can be applied with $k=5$  [see \cite{yang2010}]. Consequently, there exists a design $\xi^{+}$ with $3M$ support points and each component $\xi^{+}_i$ contains  $d^{(i)}_{\max}$ $i=1, \ldots, M$.

\bigskip
\noindent
{\bf Proof of Theorem \ref{theo:doptplac}:}
For the sake of transparency we restrict ourselves to the case $M = 2$ such that $m= Mq+1= 2q+1$.
We use  the equivalence Theorem \ref{theo:aequiv} to establish the $D$-optimality of the design $\xi^{\star}$.  In the present situation this means
that the $D$-optimality of the design $\xi^{\star}$  defined in \eqref{eq:Doptplac} for model  \eqref{eq:loc} with assumption \eqref{eq:placassump} can be proved by
 checking  the two inequalities
  \bea
 \kappa_1(t,\xi^{\star}, \theta) &=& \tfrac{1}{\sigma^2_1} \big(1, \eta^T_0(t, \theta^{(1)}_2), 0^T_q   \big)M^{-1}(\xi^{\star}, \theta)   \big(1, \eta^T_0(t, \theta^{(1)}_2), 0^T_q   \big)^T  \leq 2q+1, \, \,  t \in [0, d^{(1)}_{max}] ,~~~\label{eq:kappam1}\\
 \kappa_2 (t,\xi^{\star}, \theta) &= & \tfrac{1}{\sigma^2_2} \big(1, 0^T_q, \eta^T_0(t, \theta^{(2)}_2)   \big) M^{-1}(\xi^{\star}, \theta)  \big(1, 0^T_q, \eta^T_0(t, \theta^{(2)}_2)   \big)^T \leq 2q+1,  \,\, t \in [0, d^{(2)}_{max}],~~~
 \label{eq:kappam2}
 \eea
 where $\eta_0(d, \theta^{(i)}_2) = \tfrac{\partial}{\partial \theta^{(i)}_2} f_0 (d, \theta^{(i)}_2)$. 
A straightforward calculation shows  that the information
of the design $\xi^{\star}$ can be represented as
 \be \label{eq:infombest}
 M(\xi^{\star}, \theta) = \frac{1}{m} X (\sigma_1, \theta^{(1)}_2, \sigma_2, \theta^{(2)}_2)  X^T(\sigma_1, \theta^{(1)}_2, \sigma_2, \theta^{(2)}_2)~,
 \ee
where the matrix $X(\sigma_1, \theta^{(1)}_2, \sigma_2, \theta^{(2)}_2)$ is given by
 \beo 
 \label{xmat}
X(\sigma_1, \theta^{(1)}_2, \sigma_2, \theta^{(2)}_2)= \begin{pmatrix} X_{11} \big ({\sigma_1} , 0,
d^{(1)}_1, \ldots, d^{(1)}_q , \theta^{(1)}_2 \big) &  X_{12} ({\sigma_2} ) \\
 0 &   X_{22} \big( {\sigma_2},d^{(2)}_1, \ldots, d^{(2)}_q, \theta^{(2)}_2 \big)
\end{pmatrix},
\eeo
and   the matrices $X_{11} \big ({\sigma} , 0,
d^{(1)}_1, \ldots, d^{(1)}_q , \theta^{(1)}_2 \big)\in \R^{(q+1) \times (q+1)} $,
$X_{22} (\sigma, d^{(2)}_1, \ldots, d^{(2)}_q, \theta^{(2)}_2 ) \in \R^{q\times q}$ and
 $X_{12} ( \sigma) \in \R^{(q+1)\times q} $  are defined by
 \begin{eqnarray*}
 \label{defblockmat1}
X_{22} (\sigma, d^{(2)}_1, \ldots, d^{(2)}_q, \theta^{(2)}_2 ) &=& \frac{1}{\sigma}\big(\eta_0(d^{(2)}_1, \theta^{(2)}_2) , \ldots , \eta_0(d^{(2)}_q, \theta^{(2)}_2)\big)~,~\\
\label{defblockmat}
 X_{11} (\sigma, 0 , d^{(1)}_1, \ldots d^{(1)}_q, \theta^{(1)}_2 ) &=&
\frac{1}{\sigma}
 \begin{pmatrix}
1 & 1^T_q \\
0_q &X_{22} (\sigma, d^{(1)}_1, \ldots, d^{(1)}_q, \theta^{(1)}_2 )
\end{pmatrix}~,~ X_{12} ( \sigma) =
\frac{1}{\sigma}
  \begin{pmatrix}
1 & \ldots & 1  \\
0_q & \ldots & 0_q
 \end{pmatrix} ~.~
 \end{eqnarray*}
Consequently, the   inverse of  $ M (\xi^{\star} , \theta) $ is obtained  as
$$ M^{-1}(\xi^{\star} , \theta) =  m (X^T (\sigma_1, \theta^{(1)}_2, \sigma_2, \theta^{(2)}_2) )^{-1}X^{-1} (\sigma_1, \theta^{(1)}_2, \sigma_2, \theta^{(2)}_2),$$
 where
 {\footnotesize
   \bea \nonumber
&& X^{-1} (\sigma_1, \theta^{(1)}_2, \sigma_2, \theta^{(2)}_2) =   \\
&& ~~~~~~~~~~~~~~~ \begin{pmatrix}
X^{-1}_{11} (\sigma, 0 , d^{(1)}_1, \ldots d^{(1)}_q, \theta^{(1)}_2 ) & -  X^{-1}_{11} (\sigma, 0 , d^{(1)}_1, \ldots d^{(1)}_q, \theta^{(1)}_2 )
X_{12}  ({\sigma_2} ) \, X^{-1}_{22} (\sigma, d^{(2)}_1, \ldots, d^{(2)}_q, \theta^{(2)}_2, \theta^{(2)}_2  )   \\
0  &  X^{-1}_{22} (\sigma, d^{(2)}_1, \ldots, d^{(2)}_q, \theta^{(2)}_2 )  \end{pmatrix}. \nonumber
 \eea}
Using these block structures the function $\kappa_1(t,  \xi^{\star}  , \theta)$ defined in  \eqref{eq:kappam1}  reduces
for  the design
$\xi^{\star}  =(\xi^{\star}_1,\xi^{\star}_2,\mu^{\star})$
to
 \beao
  \kappa_1(t, \xi^{\star} , \theta)
 &=& \tfrac{m}{\sigma^2_1} \big(1,    \eta^T_0(t, \theta^{(1)}_2) \big)
  \big ( X_{11}^{-1} (\sigma, 0 , d^{(1)}_1, \ldots d^{(1)}_q)  \big )^T    \,  X_{11}^{-1} (\sigma, 0 , d^{(1)}_1, \ldots d^{(1)}_q)
 \big(1, \eta^T_0(t, \theta^{(1)}_2) \big)^T \\
 &=& \tfrac{m}{(q+1) \sigma^2_1} \big(1,    \eta^T_0(t, \theta^{(1)}_2) \big)
M_1^{-1}(\xi_1^{\star} ,\theta^{(1)} ) \big(1,    \eta^T_0(t, \theta^{(1)}_2) \big)^T~,
 \eeao
 where $ M_1(\xi_1^{\star} ,\theta_1, \theta^{(1)}_2 ) =   \frac{1}{ \sigma^2_1}  \int_0^1\big(1,    \eta^T_0(t, \theta^{(1)}_2) \big)^T
  \big(1,    \eta^T_0(t, \theta^{(1)}_2) \big) d\xi_1^{\star}(t)$  denotes the information matrix of the design $\xi_1^{ \star}$ in the single model with parameter
 $(\theta_1, \theta_2^{(1)})$.
 Consequently, the function   $\kappa_1(t, \xi^{\star}   , \theta_1, \theta^{(1)}_2)$  only depends on the first component $\xi_1^\star$
 and  is  proportional to  the left-hand side
of  the standard equivalence theorem for $D$-optimality for  the single model.
The inequality $  \kappa_1(t,  \xi^{\star}  , \theta) \le m $  for all $t\in [0, d^{(1)}_{\max}]$ follows from the fact that the design $\xi^{\star}_1$ given in \eqref{eq:doptsingle} is locally $D$-optimal for the single model with parameter $(\theta_1, \theta_2^{(1)})$ and  this proves \eqref{eq:kappam1}.  \\
 In order to show  the remaining inequality \eqref{eq:kappam2}  for all $t \in [0, d^{(2)}_{\max}]$ we use the fact that the information matrix in \eqref{eq:infombest} can be represented as
 \beo
 M( \xi^{\star}  , \theta) = S
 X ( \sigma_2, \theta^{(2)}_2,\sigma_1, \theta^{(1)}_2)
 \diag( \tfrac{\sigma^2_2}{m \sigma^2_1}, \tfrac{1}{m},  \ldots, \tfrac{1}{m} )
 X^T ( \sigma_2, \theta^{(2)}_2,\sigma_1, \theta^{(1)}_2) S,
 \eeo
  where $S$ denotes  a $m \times m$ permutation matrix, defined by
$$
S=\begin{pmatrix} 1 & 0^T_q & 0^T_q  \\  0_q & 0_{q\times q} & I_{q\times q}  \\ 0_q & I_{q\times q} & 0_{q \times q}  \end{pmatrix} , 
$$
 $0_{q\times q}$ denotes a matrix with all entries equal to zero and  $I_{q\times q}$ the $q \times q$ identity matrix.
Observing that
$S h_2(t) =   \tfrac{1}{\sigma_2}\bigl(1,  \eta^T_0(t, \theta^{(2)}_2), 0^T_q\bigr)^T$
it follows that
the  function $\kappa_2(t, \xi^{\star}  , \theta)$  in  \eqref{eq:kappam2}  can be represented as
\beao
\kappa_2(t, \xi^{\star}  , \theta)
&=&\tfrac{1}{\sigma^2_2} \big(1,    \eta^T_0(t, \theta^{(2)}_2) \big)
\big (X_{11}^{-1} \big ({\sigma_2} , 0,
 d^{(2)}_1, \ldots, d^{(2)}_q, \theta^{(2)}_2   \big)\big)^T
\Bigl[m I_{q+1 \times q+1}  \\
&& \qquad  ~~-~m(1-\tfrac{\sigma^2_1}{\sigma^2_2}) \diag( 1 , 0_q) \Bigr] X_{11}^{-1} \big ({\sigma_2} , 0, d^{(2)}_1, \ldots, d^{(2)}_q ,\theta^{(2)}_2  \big)
\big(1,    \eta^T_0(t, \theta^{(2)}_2) \big)^T \\
&=&\tfrac{m}{(q+1) \sigma^2_2} \big(1,    \eta^T_0(t, \theta^{(2)}_2) \big) M_2^{-1}(\tilde\xi_2, \theta_1,  \theta^{(2)}_2)  \big(1,    \eta^T_0(t, \theta^{(2)}_2) \big)^T \\
&&- m(1-\tfrac{\sigma^2_1}{\sigma^2_2})  \tfrac{1}{\sigma^2_2}\Bigl[ ( 1, 0^T_q) X_{11}^{-1} \big ({\sigma_2} ,  0, d^{(2)}_1, \ldots, d^{(2)}_q , \theta^{(2)}_2  \big)  \big(1,    \eta^T_0(t, \theta^{(2)}_2) \big)^T \Bigr]^2,
\eeao
 where  $M_2(\tilde\xi_2, \theta_1, \theta^{(2)}_2)  $ is the information matrix of the design  $\tilde \xi_2$ given by \eqref{eq:doptsingle} for the single model.
The  first term of this expression is proportional to the left hand side of the equivalence  theorem corresponding to  the $D$-optimality in the single model with parameter $(\theta_1, \theta_2^{(2)})$.
 Moreover, it follows that  the   design  $\tilde \xi_2 $  is  $D$-optimal for the single model with parameter $(\theta_1, \theta_2^{(2)})$, which implies that the first term is always smaller than $m$.
By  the assumption $\sigma^2_1 \leq \sigma^2_2$ we obtain that the second term of this expression is nonpositive,
which shows  $\kappa_2(t,  \xi^{\star}  , \theta)\leq m$ for all $t \in [0, d^{(2)}_{\max}]$. This proves the inequality
 \eqref{eq:kappam2}  and completes the proof of Theorem  \ref{theo:doptplac} in the case $M=2$.

\bigskip \noindent
{\bf Proof of Corollary \ref{theo:best5point}:}
The locally $D$-optimal designs for the (single) Emax, the linear-in-log and the exponential model were calculated by \cite{detkis2010}. The corollary now follows by an application of Theorem \ref{theo:doptplac}.

\bigskip\noindent
 {\bf Proof of Lemma \ref{theo:suppoints}:} 
 Let $\xi^\star= (\xi^\star_1, \xi^\star_2, \mu^\star)$ denote the locally $D$-optimal design
for the Emax, the linear-in-log or the exponential model. Since the information matrix $M(\xi^{\star}, \theta)$ of a locally $D$-optimal design  must be nonsingular
one can easily deduce the following implications
\begin{align}
 &|\supp(\xi^\star_1)| + |\supp(\xi^\star_2) | \geq 4   \label{eq:cond1}\\
& \mbox{If } |\supp(\xi^\star_1)| + |\supp(\xi^\star_2) | = 4 \,,  ~ \mbox{ then }  0 \notin \supp(\xi^\star_1) \cap \supp(\xi^\star_2)  \label{eq:cond2}\\
&  \mbox{If }  |\supp(\xi^\star_i)| = 1 \, ,~ \mbox{ then }     0 \notin \supp(\xi^\star_i),  \, i = 1, 2.  \nonumber 
\end{align}
Moreover, it follows by Corollary \ref{kor:garza} that the locally $D$-optimal design has at most $5$ support points for the Emax and the linear-in-log model and at most $6$ support points for the exponential model. This proves Assertion (A1) and (B).
Assertion (A2) also follows by Corollary \ref{kor:garza}. \\
For a proof of (A3) we note that $(|\supp(\xi^\star_1)|, |\supp(\xi^\star_2)|) \in \{(1, 3), (2, 2), (3, 1) \}$   if the locally $D$-optimal design is given by a  design in $\Xi^4_2$. 
If $(|\supp(\xi^\star_1)|, |\supp(\xi^\star_2)|) = (1, 3)$, $\xi^{\star}_2$ must contain the boundary points $0, d^{(2)}_{\max}$, otherwise it could be improved with respect to the Loewner ordering (see Theorem \ref{theo:garza}).
If $(|\supp(\xi^\star_1)|, |\supp(\xi^\star_2)|) = (2, 2)$, both designs must contain at least one of the boundary points, otherwise $I(\xi^{\star}_i) = 2$ ($i=1,2$) and the designs could be improved with respect to the Loewner ordering (see again Theorem \ref{theo:garza}). Using \eqref{eq:cond2} it follows that at least one of the designs contains the corresponding upper boundary point. If $(|\supp(\xi^\star_1)|, |\supp(\xi^\star_2)|) = (3, 1)$, $\xi^{\star}_1$ must contain the boundary points $0, d^{(1)}_{\max}$, otherwise it could be improved with respect to the Loewner ordering (see Theorem \ref{theo:garza}).
Assertion (A3) now follows. 

\bigskip \noindent
{\bf Proof of Theorem \ref{theo:best4point}: } For the sake of brevity we restrict the discussion to the Emax model.
The proof consists of two  steps. At first we show that it is sufficient to prove the result on the design space $[0,1]$.
Secondly, we determine the $D$-optimal  design in the class $\Xi_2^4$.\\
\noindent
{\bf (1) } Recall the definition of the information matrix in \eqref{eq:info} (with $M=2$) in
 model \eqref{eq:Emax2a} with the parameter vector $\theta= (\vartheta_{11}, \vartheta_{12}, \theta^{(1)}_2, \theta^{(2)}_2)^T \in \R^4$.
Let $\xi= (\xi_1, \xi_2, \mu)$ denote an arbitrary design with  components  $\xi_1$ and $\xi_2$ defined on the design space $[0,d_{\max}^{(1)}]$ and $[0,d_{\max}^{(2)}]$,
respectively, and denote by  $\tilde \xi_1$ and $\tilde \xi_2$  the corresponding measures on the interval $[0,1]$  induced by the transformation
$ t \to t/d^{(i)}_{\max}$ ($i=1,2$).  Now a straightforward calculation gives
 \beao
 M(\xi, \theta) &=& \lambda \int^{d_{\max}^{(1)}}_0  h_1(t_1) h^T_1(t_1) d\xi_1(t_1) + (1-\lambda) \int^{d_{\max}^{(2)}}_0  h_2(t_2) h^T_2(t_2) d\xi_2(t_2),  \\
&=& \lambda \int_0^1 \tilde h_1(t_1)  \tilde h^T_1(t_1) d \tilde \xi_1(t_1) + (1-\lambda) \int_0^1  \tilde h_2(t_2)  \tilde h^T_2(t_2) d \tilde  \xi_2(t_2) ,
 \eeao
where
 \beao
\tilde  h^T_1(t_1) &=& \frac{1}{\sigma_1} \Bigr(1, \frac{t_1}{t_1 + {\theta^{(1)}_2}/{d^{(1)}_{\max}}},  \frac{-t_1}{(t_1 + {\theta^{(1)}_2}/{d^{(1)}_{\max}})^2}, 0\Bigl) P , \\
 \tilde h^T_2(t_2)  &=& \frac{1}{\sigma_2}\Bigl(1, \frac{t_2}{t_2 +{\theta^{(2)}_2}/{d^{(2)}_{\max}}},  0, \frac{-t_2}{(t_2 + {\theta^{(2)}_2}/{d^{(2)}_{\max}})^2} \Bigr) P~,
 \eeao
 where  $$P ={\rm diag} \big(1,1, \tfrac{\vartheta_{1}}{d^{(1)}_{\max}} , \tfrac{\vartheta_{2}}{d^{(2)}_{\max}} \big)
$$
is a $4 \times 4$ matrix.
This shows that the  components $\xi_1^\star$  and $\xi_2 ^\star$ of the  locally $D$-optimal design
for the Emax model on the  design spaces  $[0,d_{\max}^{( 1)}] $ and  $[0,d_{\max}^{(2)}] $   can be obtained
 by a linear transformation of the corresponding locally $D$-optimal designs
on the design space $[0,1]$, where the parameters in the Emax model are given by
$\tilde \theta = \bigl(\theta_1,1, {\theta^{(1)}_2}/{d^{(1)}_{\max}} ,  1, {\theta^{(2)}_2}/{d^{(2)}_{\max}} \bigr)^T$.
Therefore it is sufficient to consider the case  ${\cal X}_1= {\cal X}_2 = [0,1]$ in the following discussion.\\
{\bf (2)} According to {\bf (1)} we restrict ourselves to the case $d_{\max}^{(1)} =  d_{\max}^{(2)} =1$  and  $\theta_1 =1$.
 Therefore the main assumption of the theorem reduces to $0 < \theta^{(1)}_2 < \theta^{(2)}_2 < 1$.
Because a $D$-optimal design in the class $\Xi^2_4$ must have a nonsingular information matrix it
follows that
\beo
|\supp(\xi_i)| \geq 1, \quad i=1, 2~.
\eeo
By an application of Lemma   \ref{lem:doptweights}  we obtain the following candidates for the $D$-optimal  design in the class $\Xi^2_4$
\begin{align}
 \xi^a_1 &= \begin{pmatrix} d^{(1)}_1 & d^{(1)}_2 & d^{(1)}_3 \\ \tfrac{1}{3} &   \tfrac{1}{3}  &  \tfrac{1}{3}  \end{pmatrix}, \quad \xi^a_2= \begin{pmatrix} d^{(2)}_1  \\ 1 \end{pmatrix}, \quad \mu^a = \begin{pmatrix} 1& 2 \\ \tfrac{3}{4} & \tfrac{1}{4}  \end{pmatrix} \label{eq:canda4} ~, \\
  \xi^b_1 &= \begin{pmatrix} d^{(1)}_1 & d^{(1)}_2 \\ \tfrac{1}{2} & \tfrac{1}{2} \end{pmatrix}, \quad \xi^b_2= \begin{pmatrix} d^{(2)}_1 &d^{(2)}_2 \\\tfrac{1}{2} & \tfrac{1}{2}  \end{pmatrix}, \quad \mu^b= \begin{pmatrix} 1& 2 \\ \tfrac{1}{2} & \tfrac{1}{2}  \end{pmatrix} \label{eq:candc4}, \\
  \xi^c_1 &= \begin{pmatrix} d^{(1)}_1 \\1  \end{pmatrix}, \quad \xi^c_2= \begin{pmatrix} d^{(2)}_1 &d^{(2)}_2 &d^{(2)}_3\\\tfrac{1}{3} &   \tfrac{1}{3}  &  \tfrac{1}{3}  \end{pmatrix}, \quad \mu^c = \begin{pmatrix} 1& 2 \\ \tfrac{1}{4} &  \tfrac{3}{4} \end{pmatrix}  \label{eq:candb4}.
\end{align}
Next we evaluate the determinants for these three candidate designs and maximize them with respect to the support points. For
example we obtain for the design $D$-optimal $\xi^a =( \xi^a_1,  \xi^a_2, \mu^a)$
\beao
\det(M(\xi^a, \theta)) &= &
\bigl(\tfrac{1}{\sigma_1^2}\bigl)^3 \bigl( \tfrac{1}{\sigma^2_2} \bigr) \bigl( \tfrac{1}{4} \bigr)^4  \mbox{det}^2 \begin{pmatrix}
 1 & 1 &1 &1 \\
\frac{d^{(1)}_1}{d^{(1)}_1 + \theta^{(1)}_2} & \frac{d^{(1)}_2}{d^{(1)}_2+ \theta^{(1)}_2} &\frac{d^{(1)}_3}{d^{(1)}_3+ \theta^{(1)}_2} &\frac{d^{(2)}_1}{d^{(2)}_1 + \theta^{(2)}_2} \\
\frac{-d^{(1)}_1}{(d^{(1)}_1 + \theta^{(1)}_2)^2} &\frac{-d^{(1)}_2}{(d^{(1)}_2+ \theta^{(1)}_2)^2} & \frac{-d^{(1)}_3}{(d^{(1)}_3+ \theta^{(1)}_2)^2} &0 \\
0 &0 & 0 &\frac{-d^{(2)}_1}{(d^{(2)}_1 + \theta^{(2)}_2)^2}
\end{pmatrix} \\
&&\\
&=&  \bigl(\tfrac{1}{\sigma_1^2}\bigl)^3 \bigl( \tfrac{1}{\sigma^2_2} \bigr) \bigl( \tfrac{1}{4} \bigr)^4 \bigl(\theta^{(1)}_2 \bigr)^4  \Bigl( \tfrac{(d^{(1)}_2-d^{(1)}_1)(d^{(1)}_3-d^{(1)}_1)(d^{(1)}_3-d^{(1)}_2)}{(d^{(1)}_1 + \theta^{(1)}_2 )^2(d^{(1)}_2+\theta^{(1)}_2  )^2(d^{(1)}_3+ \theta^{(1)}_2 )^2} \Bigr)^2 \Bigl(\tfrac{d^{(2)}_1}{(d^{(2)}_1 + \theta^{(2)}_2)^2} \Bigr)^2.
\eeao
Observing that $\theta^{(2)}_2 < 1$ we obtain  that the factor $\bigl(\tfrac{d^{(2)}_{1}}{(d^{(2)}_{1} + \theta^{(2)}_2)^2}\bigr)^2$  is maximized
for  $  d_1^{(2)} = \theta^{(2)}_2 $,  and it follows by a straightforward calculation that the  support points of the design
$\xi^a_1$ and $\xi_2^a$  in \eqref{eq:canda4} maximizing the determinant are given by
\beo 
 \label{eq:supp4a}
d^{(1)}_1 = 0, \quad d^{(1)}_2 = \tfrac{\theta^{(1)}_2}{1+ 2\theta^{(1)}_2} , \quad d^{(1)}_3 = 1,  \qquad  d^{(2)}_1 = \theta^{(2)}_2,
\eeo
respectively.
The resulting determinant of the information matrix of the corresponding  design, say  $\xi^{a,\star}$,  is
obtained as
\beo
 \label{eq:det4a}
\det(M(\xi^{a,\star}, \theta))= \frac{\bigl(\tfrac{1}{\sigma_1^2}\bigl)^3 \bigl( \tfrac{1}{\sigma^2_2} \bigr)^1 \bigl( \tfrac{1}{4} \bigr)^8 }{\bigl(\theta^{(1)}_2 \theta^{(2)}_2\bigr)^2 \bigl( 1+ \theta^{(1)}_2\bigr)^6}.
\eeo
Analogously, we get that the support points
\beo 
 \label{eq:supp4b}
d^{(1)}_1 =  \theta^{(1)}_2, \qquad d^{(2)}_1 = 0, \quad d^{(2)}_2= \tfrac{\theta^{(2)}_2}{1+ 2\theta^{(2)}_2} , \quad d^{(2)}_3 = 1,
\eeo
yield a maximal determinant  for the designs $\xi_1^c$, $\xi_2^c$ in \eqref{eq:candb4}, respectively,   and the
  determinant of the information matrix of the corresponding design, say  $\xi^{c,\star} $
is given by
\beo 
 \label{eq:det4b}
\det(M(\xi^{c,\star} , \theta))= \frac{ \bigl(\tfrac{1}{\sigma_1^2}\bigl) \bigl( \tfrac{1}{\sigma^2_2} \bigr)^3 \bigl( \tfrac{1}{4} \bigr)^8 }{\bigl(\theta^{(1)}_2 \theta^{(2)}_2\bigr)^2 \bigl( 1+ \theta^{(2)}_2\bigr)^6}.
\eeo
Finally, we consider the determinant of the   candidate design \eqref{eq:candc4}, that is
\beo 
 \label{eq:Emax4c}
\det(M(\xi^b, \theta) ) =  \bigl(\tfrac{1}{\sigma_1^2}\bigl)^2 \bigl( \tfrac{1}{\sigma^2_2} \bigr)^2 \bigl( \tfrac{1}{4} \bigr)^4 \frac{(d^{(1)}_2 -d^{(1)}_1)^2 (d^{(2)}_2 -d^{(2)}_1)^2 (d^{(2)}_2d^{(2)}_1 \bigl(\theta^{(1)}_2\bigr)^2 - d^{(1)}_2d^{(1)}_1 \bigl(\theta^{(2)}_2\bigr)^2 )^2}{(d^{(1)}_1 + \theta^{(1)}_2)^4(d^{(1)}_2 +\theta^{(1)}_2)^4(d^{(2)}_1 + \theta^{(2)}_2)^4(d^{(2)}_2 + \theta^{(2)}_2)^4}.
\eeo
We assume without loss of generality  that $d^{(i)}_1 < d^{(i)}_2$ $i=1, 2$.
Note that for $D$-optimality of the design $\xi^b$ the smallest support points $d^{(1)}_1$ and  $d^{(2)}_1$
of the components  $\xi^b_1$ and  $\xi^b_2$
must satisfy  $d^{(1)}_1 +d^{(2)}_1 > 0$  (otherwise the determinant vanishes).  Consequently
there exist two possible cases for the design $\xi^{b,\star} $ corresponding to the cases $d^{(1)}_1 =0$ or $d^{(2)}_1=0$, namely
\begin{align*} 
 \xi^{b_1}_1 &= \begin{pmatrix} 0 & d^{(1)}_2 \\ \tfrac{1}{2} & \tfrac{1}{2} \end{pmatrix}, \quad \xi^{b_1}_2= \begin{pmatrix} d^{(2)}_1 &d^{(2)}_2 \\ \tfrac{1}{2} & \tfrac{1}{2} \end{pmatrix}, \quad \mu^{b_1}= \begin{pmatrix} 1& 2 \\ \tfrac{1}{2} & \tfrac{1}{2} \end{pmatrix}, \\ 
  \xi^{b_2}_1 &= \begin{pmatrix} d^{(1)}_1 & d^{(1)}_2 \\ \tfrac{1}{2} & \tfrac{1}{2} \end{pmatrix}, \quad \xi^{b_2}_2= \begin{pmatrix} 0 &d^{(2)}_2 \\ \tfrac{1}{2} & \tfrac{1}{2} \end{pmatrix}, \quad \mu^{b_2}= \begin{pmatrix} 1& 2 \\ \tfrac{1}{2} & \tfrac{1}{2}\end{pmatrix}. 
\end{align*}
Now a straightforward calculation gives for the design $ \xi^{b_1}_1$
\beao
\det(M(\xi^{b_1}, \theta)) &=&
\bigl(\tfrac{1}{\sigma_1^2}\bigl)^2 \bigl( \tfrac{1}{\sigma^2_2} \bigr)^2 \bigl( \tfrac{1}{4} \bigr)^4  \mbox{det}^2 \begin{pmatrix}
 1 & 1 &1 &1 \\
0 & \frac{d^{(1)}_{2}}{d^{(1)}_{2} + \theta^{(1)}_2} &\frac{d^{(2)}_{1}}{d^{(2)}_{1} + \theta^{(2)}_2} & \frac{d^{(2)}_{2}}{d^{(2)}_{2} + \theta^{(2)}_2} \\
0 &\frac{-d^{(1)}_{2}}{(d^{(1)}_{2} + \theta^{(1)}_2)^2} & 0 &0 \\
0 &0 &  \frac{-d^{(2)}_{1}}{(d^{(2)}_{1} + \theta^{(2)}_2)^2} &\frac{-d^{(2)}_{2}}{(d^{(2)}_{2} + \theta^{(2)}_2)^2}
\end{pmatrix} \\
&&\\
&=& \bigl(\tfrac{1}{\sigma_1^2}\bigl)^2 \bigl( \tfrac{1}{\sigma^2_2} \bigr)^2 \bigl( \tfrac{1}{4} \bigr)^4 \Bigl(\tfrac{(d^{(2)}_2-d^{(2)}_1)d^{(2)}_1 d^{(2)}_2}{(d^{(2)}_1 + \theta^{(2)}_2 )(d^{(2)}_2+\theta^{(2)}_2  )} \Bigr)^2 \Bigl(\tfrac{d^{(1)}_2}{(d^{(1)}_2 + \theta^{(1)}_2)^2} \Bigr)^2
\eeao
Maximizing with respect to the support points yields
\beo
d^{(1)}_1 = 0, \quad d^{(1)}_2 = \theta^{(1)}_2, \qquad d^{(2)}_1  = \tfrac{\theta^{(2)}_2}{1+ 2\theta^{(2)}_2}, \quad d^{(2)}_2 = 1,
\eeo
and we obtain for the determinant of the design $\xi^{b_1,\star} $
\bea \label{eq:det4c1}
\det(M(\xi^{b_1,\star}, \theta)) &=&  \frac{ \bigl(\tfrac{1}{\sigma_1^2}\bigl)^2 \bigl( \tfrac{1}{\sigma^2_2} \bigr)^2 \bigl( \tfrac{1}{4} \bigr)^8}{\bigl(\theta^{(1)}_2 \theta^{(2)}_2\bigr)^2 \bigl( 1+ \theta^{(2)}_2\bigr)^6}
\eea
A similar optimization of  the determinant of the information matrix of the design $\xi^{b_2,\star}$ gives
\bea
\det(M(\xi^{b_2,\star}, \theta)) &=& \frac{\bigl(\tfrac{1}{\sigma_1^2}\bigl)^2 \bigl( \tfrac{1}{\sigma^2_2} \bigr)^2 \bigl( \tfrac{1}{4} \bigr)^8}{\bigl(\theta^{(1)}_2 \theta^{(2)}_2\bigr)^2 \bigl( 1+ \theta^{(1)}_2\bigr)^6} ~. \label{eq:det4c}
\eea
By a comparison of \eqref{eq:det4c1} and \eqref{eq:det4c}
 it follows that the determinant of  design $\xi^{b_2,\star}$ is always larger than the determinant of  the design $\xi^{b_1,\star}$,
 since $\theta^{(1)}_2 < \theta^{(2)}_2$. Finally,  the assertion of the theorem follows by straightforward calculations
 comparing the determinants of the designs  $\xi^{a,\star}, \xi^{b,\star}, \xi^{c_2,\star}$ in the different scenarios for the ratio
 $r= \tfrac{\sigma_1^2}{\sigma_2^2}$.

\medskip

\noindent
{\bf Proof of Theorem   \ref{theo:Emaxdopt}:}
By similar arguments as given in the proof of Theorem
\ref{theo:best4point} we obtain that it is sufficient to consider the case $d_{\max}^{(1)}= d_{\max}^{(2)}=1$.\\
\noindent
{\bf (1)} In the  case    $r \leq 1$   it follows  from Theorem \ref{theo:aequiv} that the
design $\xi^{a, \star}$ defined  in  \eqref{eq:best4pointa} is locally $D$-optimal for model \eqref{eq:Emax2a}
if and only if  the two inequalities
 \bea
 \kappa_1(t, \xi^{a, \star}, \theta) &=& \tfrac{1}{\sigma^2_1} \Bigl(1, \tfrac{t}{t+ \theta^{(1)}_2} , \tfrac{-t}{(t+ \theta^{(1)}_2)^2}, 0   \Bigr)M^{-1}(\xi^{a, \star}, \theta)  \Bigl(1, \tfrac{t}{t+ \theta^{(1)}_2} , \tfrac{-t}{(t+ \theta^{(1)}_2)^2}, 0  \Bigr)^T  \leq 4 \label{eq:kappa41}\\
 \kappa_2 (t, \xi^{a, \star}, \theta) &= & \tfrac{1}{\sigma^2_2} \Bigl(1,\tfrac{t}{t+ \theta^{(2)}_2} ,0 , \tfrac{-t}{(t+ \theta^{(2)}_2)^2}  \Bigr)M^{-1}(\xi^{a, \star}, \theta)  \Bigl(1,\tfrac{t}{t+ \theta^{(2)}_2} , 0, \tfrac{-t}{(t+ \theta^{(2)}_2)^2}  \Bigr)^T \leq 4 \label{eq:kappa42}
 \eea
hold for all $t\in [0,1]$ [see  Theorem \ref{theo:aequiv}]. The information matrix of the design  $\xi^{a, \star}$ can be represented as
 \beo 
  \label{eq:info4best}
 M(\xi^{a, \star}, \theta) =  \tfrac{1}{4} \tilde X (\sigma_1,\theta_2^{(1)},\sigma_2, \theta_2^{(2)})
  \tilde X ^T(\sigma_1,\theta_2^{(1)},\sigma_2, \theta_2^{(2)}) ,
 \eeo
 where
$$
 \tilde X (\sigma_1,\theta_2^{(1)},\sigma_2, \theta_2^{(2)}) =
 \begin{pmatrix} X_{11}(\sigma_1, 0, \tfrac{\theta^{(1)}_2}{2\theta^{(1)}_2 +1}, 1)  & \tilde
 X_{12}(\sigma_2, \theta^{(2)}_2) \\ 0 &  \tilde X_{22}(\sigma_2, \theta^{(2)}_2) \end{pmatrix}
 $$
and the matrices $X_{11}$, $\tilde X_{12} $ and $\tilde X_{22} $ are defined by
\beao
X_{11}(\sigma_1,d^{(1)}_1, d^{(1)}_2, d^{(1)}_3)&=& \begin{pmatrix}
 1 & 1 &1  \\
\frac{d^{(1)}_1}{d^{(1)}_1 + \theta^{(1)}_2} & \frac{d^{(1)}_2}{d^{(1)}_2+ \theta^{(1)}_2} &\frac{d^{(1)}_3}{d^{(1)}_3+ \theta^{(1)}_2} \\
\frac{-d^{(1)}_1}{(d^{(1)}_1 + \theta^{(1)}_2)^2} &\frac{-d^{(1)}_2}{(d^{(1)}_2+ \theta^{(1)}_2)^2} & \frac{-d^{(1)}_3}{(d^{(1)}_3+ \theta^{(1)}_2)^2}  \end{pmatrix} ~, \\
\tilde X_{12}(\sigma_2, \theta^{(2)}_2) &=&
\begin{pmatrix}  \tfrac{1}{\sigma_2 }\\
 \tfrac{1}{2\sigma_2}
 \end{pmatrix}~,~  \tilde X_{12}(\sigma_2, \theta^{(2)}_2) =
  \begin{pmatrix}
 \tfrac{-1}{4\theta^{(2)}_2\sigma_2 }
  \end{pmatrix}~,
\eeao

 respectively.  A straightforward calculation of the inverse of the matrix $\tilde X $ yields
\beo
\tilde X^{-1}(\sigma_1, \theta^{(1)}_2, \sigma_2, \theta^{(2)}_2) = \begin{pmatrix}  X^{-1}_{11}(\sigma_1, 0, \tfrac{\theta^{(1)}_2}{2\theta^{(1)}_2 +1}, 1)  & - X^{-1}_{11}(\sigma_1, 0, \tfrac{\theta^{(1)}_2}{2\theta^{(1)}_2 +1}, 1) \tilde X_{12}(\sigma_2, \theta^{(2)}_2) \tilde
X^{-1}_{22}(\sigma_2, \theta^{(2)}_2)\\ 0 &  \tilde X^{-1}_{22}(\sigma_2, \theta^{(2)}_2)\end{pmatrix},
 \eeo
 and we obtain for the function
 $\kappa_1(t, \xi^{a,\star}, \theta)$ in \eqref{eq:kappa41}
  the representation
 \beao
  \kappa_1(t, \xi^{a, \star}, \theta)
 &=& \tfrac{4}{3 \sigma^2_1} \Bigl(1, \tfrac{t}{t+ \theta^{(1)}_2} , \tfrac{-t}{(t+ \theta^{(1)}_2)^2}   \Bigr)
 3X^{-T}_{11}(\sigma_1, 0, \tfrac{\theta^{(1)}_2}{2\theta^{(1)}_2 +1}, 1) X^{-1}_{11}(\sigma_1, 0, \tfrac{\theta^{(1)}_2}{2\theta^{(1)}_2 +1}, 1)
 \Bigl(1, \tfrac{t}{t+ \theta^{(1)}_2} , \tfrac{-t}{(t+ \theta^{(1)}_2)^2} \Bigr)^T \\
 &=& \tfrac{4}{3 \sigma^2_1} \Bigl(1, \tfrac{t}{t+ \theta^{(1)}_2} , \tfrac{-t}{(t+ \theta^{(1)}_2)^2}   \Bigr)
M_1^{-1}(\xi^{a, \star}_1,\theta^{(1)} ) \Bigl(1, \tfrac{t}{t+ \theta^{(1)}_2} , \tfrac{-t}{(t+ \theta^{(1)}_2)^2} \Bigr)^T,
 \eeao
 where $ M_1 (\xi^{a, \star}_1,\theta^{(1)} ) $  is the information matrix of the design $\xi^{a, \star}_1$
 in the Emax model with parameter vector $(\theta_1, \theta^{(1)}_2)^T$.
Because  the design $\xi^{a, \star}_1$ given in \eqref{eq:best4pointa}
is in fact locally $D$-optimal for this  model, it follows that $ \kappa_1(t, \xi^{a, \star}, \theta) \le 4$, which proves the
first inequality of the equivalence theorem.   \\
 In order to show that the inequality in \eqref{eq:kappa42} holds for all $t \in [0, 1]$ we
 note that this inequality is equivalent to
 \bea \label{pineq}
 P(t) &=& (t+ \theta^{(2)}_2)^4 \bigl(\kappa_2(t, \xi^{a, \star}, \theta) - 4 \bigr) =  \alpha_{21} t^4 + \alpha_{22} t^3 + \alpha_{23} t^2 + \alpha_{24} t + \alpha_{25}  ~\leq ~0,
 \eea
where the last identity defines the coefficients $\alpha_{2j}$ in an obvious manner. For example,  the
 leading coefficient and the intercept are given by
\beao
\alpha_{21} &= & \tfrac{1}{\sigma^2_2} (1, 1, 0, 0) M^{-1}(\xi^{a, \star},  \theta) (1, 1, 0, 0)^T - 4 =  24 r \theta^{(1)}_2 (\theta^{(1)}_2+ 1) (2\theta^{(1)}_2 +1)^2 -4 \bigl(1-   r \bigr), \\
\alpha_{25} &=& (\theta^{(2)}_2)^4 ( \tfrac{1}{\sigma^2_2} (1, 0, 0, 0) M^{-1}(\xi^{a, \star}, \theta) (1, 0, 0, 0)^T - 4)=  4(\theta^{(2)}_2)^4 \bigl(r -1\bigr)  ,
\eeao
respectively.
Consider the case $r < 1$ (the case $r \leq 1$ is finally obtained considering the
corresponding limit) and    note that $P(0)  = \alpha_{25}  < 0$. Consequently,  \eqref{pineq} holds if  either there
are no roots of $P$ in the interval  $(0,1)$  or all roots of $P$ in the interval $(0,1)$ have multiplicity $2$.
The roots of $P(t)$ are easily calculated as
 \beo
d^{(2)}_1= \theta^{(2)}_2 \quad \tilde{d}_1= \theta^{(2)}_2 \tfrac{3 + r \theta^{(1)}_2(1+ 2 \theta^{(1)}_2) - \sqrt{s(\theta^{(1)}_2)}}{\tfrac{1}{4} \alpha_{21}}\quad \tilde{d}_2=   \theta^{(2)}_2 \tfrac{3 + r \theta^{(1)}_2(1+ 2 \theta^{(1)}_2) + \sqrt{s(\theta^{(1)}_2)}}{\tfrac{1}{4}  \alpha_{21}}~,
 \eeo
 where we use the notation
 \beo
 s(\theta^{(1)}_2)= 8 - r^2 (1+ \theta^{(1)}_2)^2 (1+ 4 \theta^{(1)}_2 + 20  (\theta^{(1)}_2)^2 ) + 2 r (1+ 6\theta^{(1)}_2  +  21 (\theta^{(1)}_2)^2 +  24(\theta^{(1)}_2)^3 +12( \theta^{(1)}_2)^4 )~.
 \eeo
Note that $s(\theta^{(1)}_2)$ is  positive (because  $\theta^{(1)}_2 >0 $ and $r \leq 1$) and that $\theta^{(2)}_2 \in (0,1) $ is a root
 of multiplicity $2$. Moreover, $P(-\theta^{(2)}_2)>0$ (since $M^{-1}(\xi^{a,\star}, \theta)$ is positive definite),
and  it follows from  $P(0) < 0$  that $P$ has  a root in the interval  $(-\theta^{(2)}_2,0)$. This
is either $\tilde{d}_1$ or $\tilde{d}_2$ depending on the sign of the leading coefficient $ \alpha_{21}$. The inequality \eqref{pineq} holds,
if the other root  is neither in $(0,1)$. \\
 In order to check the location of the roots  $\tilde d_1$ and $\tilde d_2$ we consider the condition  \eqref{eq:Emaxdopt4a1} and the case that the right hand side of   \eqref{eq:Emaxdopt4a1} is positive. This implies that  the leading coefficient $\alpha_{21}$ is positive and the root $\tilde d_2 $ is also positive.
We obtain from the condition $\tilde d_1 \in   (-\theta_2^{(2)} ,0)$ the inequality
 \beo
3 + r \, \theta^{(1)}_2(1+ 2 \theta^{(1)}_2) < \sqrt{ s(\theta^{(1)}_2)}~.
\eeo
This gives for the second root
\beo
\tilde d_2 >  \theta^{(2)}_2 \frac{6 + 2 r\, \theta^{(1)}_2(1+ 2 \theta^{(1)}_2)}{\tfrac{1}{4}  \alpha_{21}}.
\eeo
Therefore it follows from  \eqref{eq:Emaxdopt4a1} (with positive right hand side) that  the inequality $\tilde d_2 \geq 1$ is satisfied.\\
On the other hand,  if the right hand side of \eqref{eq:Emaxdopt4a1} is negative,
the leading coefficient $\alpha_{21}$ is negative and  the conditions $P(0) < 0$ and $P(-\theta^{(2)}_2) > 0$ imply that
both roots $\tilde d_1$ and $\tilde d_2$ must be  negative, because otherwise the polynomial  $P$ does not satisfy \eqref{pineq}.
Observing that $\tilde d_2 < \tilde d_1$ in this case,  it is easy to see that  the condition \eqref{eq:Emaxdopt4a1} (with negative right hand side) implies $ \tilde d_1 < 0$. \\
Summarizing,  in the case $r \leq 1$ the  inequality \eqref{eq:Emaxdopt4a1} implies \eqref{pineq} for
all $t\in [0,1]$ and the $D$-optimality of the designs $\xi^{a,\star}$ follows by an application of Theorem \ref{theo:aequiv}.

\noindent
{\bf (2)} At first, we show that the condition  \eqref{eq:Emaxdopt4b} and $r>1$ imply that $1<  r \leq \tfrac{(1+ \theta^{(2)}_2)^6}{(1+ \theta^{(1)}_2)^6}$.  The last inequality is equivalent to $\theta^{(2)}_2\geq r^{1/6}(1+ \theta^{(1)}_2)-1$ and we have to show that:
\be \label{eq:ineqr}
\frac{   (\theta^{(1)}_2)^2 (1+ 2 \theta^{(1)}_2)^2 +r (1+ \theta^{(1)}_2)^2 (1+ 4 \theta^{(1)}_2 + 20 (\theta^{(1)}_2)^2) -1}{6 + 2 \theta^{(1)}_2(1+ 2\theta^{(1)}_2)}>r^{1/6}(1+ \theta^{(1)}_2)-1.
\ee
This inequality can be rewritten by
{\small $$(20r + 4)(\theta^{(1)}_2)^4 + (44r-\sqrt[6]{r} +1)(\theta^{(1)}_2)^3 + (29r-6\sqrt[6]{r} +5)(\theta^{(1)}_2)^2 + (6r-8 \sqrt[6]{r}+2)(\theta^{(1)}_2) + (r- 6\sqrt[6]{r} + 5) >0. $$}
Note that the coefficients of the polynomial are positive for all $r > 1$. It follows by the rule of Decartes that this polynomial has no positive roots and consequently, \eqref{eq:ineqr} is satisfied for all positive $\theta^{(1)}_2$.\\
Thus, if  $r \geq 1 $ and  the inequality \eqref{eq:Emaxdopt4b} holds, we investigate the $D$-optimality of the design  $\xi^{b, \star}$ defined  by \eqref{eq:best4pointb}
 checking  the two inequalities
 \bea
 \kappa_1(t, \xi^{b, \star}, \theta) &=& \tfrac{1}{\sigma^2_1} \Bigl(1, \tfrac{t}{t+ \theta^{(1)}_2} , \tfrac{-t}{(t+ \theta^{(1)}_2)^2}, 0   \Bigr)M^{-1}(\xi^{b, \star}, \theta)  \Bigl(1, \tfrac{t}{t+ \theta^{(1)}_2} , \tfrac{-t}{(t+ \theta^{(1)}_2)^2}, 0  \Bigr)^T  \leq 4 \label{eq:kappa41b}\\
 \kappa_2 (t, \xi^{b, \star}, \theta) &= & \tfrac{1}{\sigma^2_2} \Bigl(1, 0,\tfrac{t}{t+ \theta^{(2)}_2} ,0 , \tfrac{-t}{(t+ \theta^{(2)}_2)^2}  \Bigr)M^{-1}(\xi^{b, \star}, \theta)  \Bigl(1,\tfrac{t}{t+ \theta^{(2)}_2} , 0, \tfrac{-t}{(t+ \theta^{(2)}_2)^2}  \Bigr)^T \leq  4 \label{eq:kappa42b}
 \eea
on the interval $[0,1]$ [see  Theorem \ref{theo:aequiv}].\\
 Analogously to the proof of part (1) it can be shown  that the first inequality \eqref{eq:kappa41b} is satisfied for all $t \in [0, 1]$.
 In order to establish  the inequality  \eqref{eq:kappa42b} for all $t \in [0, 1]$ we consider the polynomial
 \beao
 P(t) &=& (t+ \theta^{(2)}_2)^4 \bigl(\kappa_2(t, \xi^{b, \star}, \theta) - 4 \bigr) =  \alpha_{21} t^4 + \alpha_{22} t^3 + \alpha_{23} t^2 + \alpha_{24} t + \alpha_{25} ,
 \eeao
where the leading coefficient and the intercept are now given by
\beao
\alpha_{21} &=& \alpha_{21}(\theta^{(1)}_2) =  4 \Bigl((\theta^{(1)}_2)^2 (1+ 2 \theta^{(1)}_2)^2 + r (1+ \theta^{(1)}_2)^2 (1+ 4 \theta^{(1)}_2 + 20 (\theta^{(1)}_2)^2) \Bigr) -4, \\
\alpha_{25} &=&  4(\theta^{(2)}_2)^4 \Bigl(\tfrac{\sigma^2_2}{\sigma^2_2} -1\Bigr) = 0.
\eeao
 Moreover, $P(-\theta^{(2)}_2)>0$ (since $M^{-1}(\xi^{b, \star}, \theta)$ is positive definite) and
  the leading coefficient $\alpha_{21}$ is always positive, since $\alpha_{21}(0 ) = 4r-4 > 0$ and $\alpha_{21}$ is increasing for $\theta^{(1)}_2 \geq 0$. The roots of $P(t)$ are given by
 \beo
d^{(2)}_1 = 0, \quad d^{(2)}_2= \theta^{(2)}_2, \quad \tilde{d}_1= \theta^{(2)}_2 \tfrac{6 + 2 \theta^{(1)}_2 ( 1 + 2 \theta^{(1)}_2)}{\tfrac{1}{4} \alpha_{21}}
\eeo
 where $d^{(2)}_2$ is a root of second order.
 Now the inequality $P(t) \le 0$ holds for all $t\in [0,1]$ if and only if  $\tilde d_1 \ge 1$.
 It is easy to see that this condition is equivalent to
\beo
 \theta^{(2)}_2 \geq \frac{   (\theta^{(1)}_2)^2 (1+ 2 \theta^{(1)}_2)^2 + r (1+ \theta^{(1)}_2)^2 (1+ 4 \theta^{(1)}_2 + 20 (\theta^{(1)}_2)^2) -1}{6 + 2 \theta^{(1)}_2(1+ 2\theta^{(1)}_2)}
\eeo
which coincides with  \eqref{eq:Emaxdopt4b}.

\noindent
{\bf (3)} At first, one can show that the condition \eqref{eq:Emaxdopt4c1} and $r>1$ imply that $r  \geq \tfrac{(1+ \theta^{(2)}_2)^6}{(1+ \theta^{(1)}_2)^6}$. Then the result 
follows by similar arguments as given in the proof of part (1), which are omitted for the sake of brevity.


\section{Appendix: Derivation of candidate models based on a
  preliminary PK/PD model} \label{sec:appendix2}
\def\theequation{8.\arabic{equation}}
\setcounter{equation}{0}

\begin{figure}[h!]
\centering
  \includegraphics[width=0.82\textwidth]{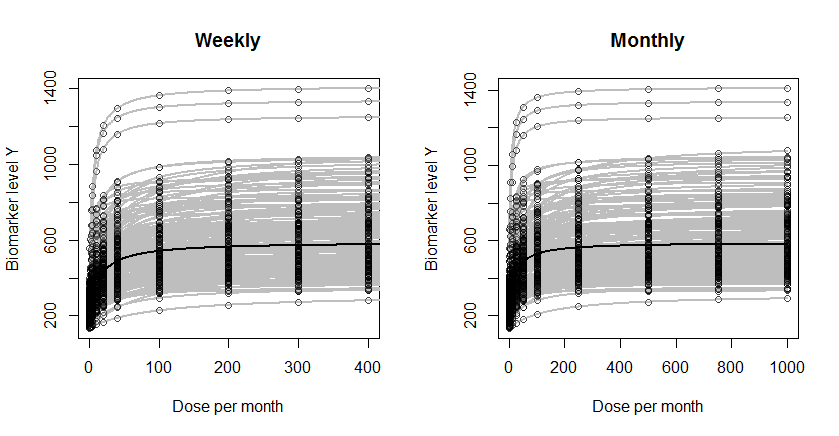}
\caption{\label{simulat}
\emph{$200$ simulated dose-response curves at timepoint T for monthly
  (group 1, right panel) and for the weekly (group 2, left panel) in grey.
  In red is the fit of the population dose-response curve.}}
\end{figure}

The PK/PD model was a nonlinear mixed effects longitudinal model
describing the PK of the drug and linking this to the PD of the drug.
The model was used to simulate longitudinal profiles per patient. The
simulation took into account parameter uncertainty from the model fit.
Then an Emax dose-group-response model was fitted to the
cross-sectional data at time $T$, that assumed that the placebo and
maximum effect of the curve are the same, but the $\mbox{ED}_{50}$ are
different in the two group. First this model was fit to the whole
population of simulations to give a population best guess, giving
$$\vartheta_{11}=  5.48, \   \vartheta_{12}=0.90, \  \theta_2^{(1)}=13.82, \
\theta_2^{(2)}= 10.46.$$
In addition $200$ individual patient profiles were simulated (see Figure \ref{simulat}) and each
individual dose-response curve was fitted at time $T$ to give $200$ parameter
estimates, representing the variability on the dose-response curve in
the population. These $200$ parameter sets
 are used to compute the distribution
of each parameter ($\vartheta_{11}$, $\vartheta_{12}$, $\theta_2^{(1)}$, $\theta_2^{(2)}$)
and their summary statistics are given in Table \ref{summary}. The logarithm of
biomarker $Y$ was modeled to achieve a better approximation through the
normal distribution.

\begin{table}[h]
\caption{\it Summary statistics}
\begin{center}
\begin{tabular}{|c|ccc|}
\hline
Parameter  & 10\% quantile & median & 90\% quantile \\
\hline
$\vartheta_{11}$ &  5.47 & 5.09 & 5.84 \\
$\vartheta_{12}$ & 0.93 & 0.66  & 1.20 \\
$\theta_2^{(1)}$ & 20.39 & 2.93 & 53.49  \\
$\theta_2^{(2)}$ & 14.99 & 2.39 & 40.40 \\
\hline
\end{tabular}
\end{center}
\end{table} \label{summary}
Based on that, we propose 4 extreme models using the 10\% and 90\%
quantile for the $\mbox{ED}_{50}$'s parameter.

\end{document}